\pdfoutput=1

\RequirePackage{fix-cm}

\documentclass[leqno, fleqn, centertags, 12pt]{article}

\usepackage{latexsym}
\usepackage{amsmath}
\usepackage{amscd}
\usepackage{amssymb}
\usepackage{verbatim}

\usepackage[T1]{fontenc}
\usepackage{latexsym}
\usepackage{manfnt}
\usepackage[sans]{dsfont}
\usepackage{palatino}
\usepackage[sc, osf]{mathpazo}
\usepackage{eulervm}

\usepackage{fullpage}

\usepackage{graphicx}

\newcommand{\nn}{\mathbold{N}}

\newcommand{\zz}{\mathbold{Z}}

\newcommand{\kk}{\Bbbk}

\newcommand{\field}{{F}}

\newcommand{\seqs}{\mathfrak{S}_{\kk}}
\newcommand{\seqq}{\mathfrak{S}}

\newcommand{\dual}{\mathcal{D}}

\newcommand{\id}{\mathop{\mathrm{id}}\nolimits}

\newcommand{\hhom}{\mathop{\mathrm{Hom}\dff}\nolimits}
\newcommand{\kker}{\mathop{Ker}\nolimits\hspace*{0.1em}}
\newcommand{\nnil}{\mathop{Nil}\nolimits\hspace*{0.1em}}
\newcommand{\nnnil}{\mathop{Nil}\nolimits\hspace*{-0.05em}}

\newcommand{\eend}{\mathop{\mbox{\textbf{End}}}\nolimits}

\newcommand{\tto}{\,{\to}\,}

\newcommand{\ssub}{\,{\subset}\,}

\newcommand{\elem}{\,{\in}\,}

\newcommand{\geqs}{\hskip.05em\,{\geq}\,}        
\newcommand{\leqs}{\,{\leq}\,\hskip.05em}        
\newcommand{\eeq}{{\hskip.05em\,{=}\,}}         
\newcommand{\nneq}{{\hskip.05em\,{\neq}\,}}      
\newcommand{\gres}{\hskip.05em\,{>}\,}

\def\hss{\hskip.025em\ }
\def\sss{\hskip.05em\ }
\def\dss{\hskip.1em\ }
\def\trs{\hskip.15em\ }
\def\qss{\hskip.2em\ }
\def\oss{\hskip.4em\ }

\def\hff{{\hskip.025em}}
\def\fff{{\hskip.05em}}
\def\dff{{\hskip.1em}}
\def\trf{{\hskip.15em}}
\def\qff{{\hskip.2em}}
\def\off{{\hskip.4em}}

\def\endss{\hspace*{0.05em}}

\def\ffdot{\hspace*{-0.1em}.\hspace*{0.4em}\ }
\def\dfdot{\hspace*{-0.2em}.\hspace*{0.4em}\ }

\def\ffcom{\hspace*{-0.1em},\hspace*{0.4em}\ }
\def\dfcom{\hspace*{-0.2em},\hspace*{0.4em}\ }

\newcommand{\nsp}{\hspace*{-0.1em}}
\newcommand{\dnsp}{\hspace*{-0.2em}}

\parskip=\the\bigskipamount

\makeatletter
\renewcommand{\@makefntext}[1]{\vspace*{0.5ex}\parindent=0em\noindent
\hspace*{-0.4em}
\hbox to 0.4em{\hss\@makefnmark}\hspace*{0.4em}{#1}
}
\makeatother

\newcounter{mysectionnumber}
\newcounter{myparnum}[mysectionnumber]
\newcommand{\mysection}[2]{\setcounter{footnote}{0}
\refstepcounter{mysectionnumber}
\section*{ \textnormal{{\themysectionnumber.} {#1}}}\label{#2}} 

\newcommand{\mynonumsection}[1]{\setcounter{footnote}{0}
\section*{\textnormal{#1}}}  

\newcommand{\mypar}[2]{\refstepcounter{myparnum}
{\vspace{\bigskipamount}\textbf{\textit{{
\themyparnum. #1}\label{#2}}}\hspace*{0.4em}}}
\renewcommand{\themyparnum}{\themysectionnumber.\arabic{myparnum}}

\newcommand{\myitpar}[1]{\vspace{\bigskipamount}\textbf{\textit{#1}}\hspace*{0.4em}}
\newcommand{\myit}[1]{\textbf{\textit{#1}}\hspace{0.0em}}
\newcommand{\mytitle}[1]{\textbf{\textit{#1}}}

\newcommand{\proof}{\vspace{-\bigskipamount}{\paragraph{{\emph{Proof}.\hspace*{0.2em}}}\hspace{-0.7em}}}

\newcommand{\eproof}{ $\blacksquare$}
\newcommand{\esubproof}{ $\square$}

\newcommand{\minus}{\hspace*{0.15em}\mbox{\rule[0.4ex]{0.4em}{0.4pt}}\hspace*{0.15em}}
\newcommand{\mminus}{\hspace*{0.15em}\mbox{\rule[0.5ex]{0.4em}{0.4pt}}\hspace*{0.15em}}
\newcommand{\plus}{\hspace*{0.15em}\mbox{\rule[0.4ex]{0.5em}{0.4pt}\hspace*{-0.27em}\rule[-0.15ex]{0.4pt}{1.2ex}\hspace*{0.27em}}\hspace*{0.15em}}
\newcommand{\pplus}{\hspace*{0.15em}\mbox{\rule[0.5ex]{0.5em}{0.4pt}\hspace*{-0.27em}\rule[-0.1ex]{0.4pt}{1.2ex}\hspace*{0.27em}}\hspace*{0.15em}}

\newcommand{\deseq}[1]{\ensuremath{D^{#1}}\hspace*{0.02em}}

\newcommand{\aang}[2]{{\langle}\dff{#1}\dff,\dff{#2}\dff{\rangle}}
\newcommand{\bico}[2]{(\dff #1\hspace*{-0.1em}\mid\hspace*{-0.1em} #2\dff)}

\begin{document}

\title{Algebra\hspace*{0.1em} of\hspace*{0.1em} 
linear\hspace*{0.1em} recurrence\hspace*{0.1em} relations\hspace*{0.1em}\\ 
in\hspace*{0.1em} arbitrary\hspace*{0.1em} characteristic} 
\date{}
\author{\textnormal{Nikolai\hspace*{0.1em} V.\hspace*{0.1em} Ivanov}}

\maketitle

\footnotetext{\hspace*{-0.5em}\copyright\ Nikolai V. Ivanov, 2014\hff,\trs 2016.\trs 
Neither the work reported in this paper{}, nor its preparation were supported 
by any governmental or non-governmental agency, foundation, or institution.}

\vspace*{30ex}

\myit{\hspace*{0em}\large Contents}\vspace*{1ex} \vspace*{\bigskipamount}\\ 
\myit{Preface}\hspace*{0.5em}  \hspace*{0.5em} \vspace*{1ex}\\
\myit{1.\phantom{*}}\hspace*{0.25em} Divided derivatives of polynomials\vspace*{0.25ex}\\
\myit{2.\phantom{*}}\hspace*{0.25em} Sequences and duality\vspace*{0.25ex}\\
\myit{3.\phantom{*}}\hspace*{0.25em} Adjoints of the left shift and of divided derivatives\vspace*{0.25ex}\\
\myit{4.\phantom{*}}\hspace*{0.25em} Endomorphisms and their eigenvalues\vspace*{0.25ex}\\
\myit{5.\phantom{*}}\hspace*{0.25em} Torsion modules and a property of free modules\vspace*{0.25ex}\\
\myit{6.\phantom{*}}\hspace*{0.25em} Polynomials and their roots\vspace*{0.25ex}\\
\myit{7.\phantom{*}}\hspace*{0.25em} The main theorems\vspace{1ex}\\
\myit{Note bibliographique}\vspace{1ex}\\
\myit{Bibliographie} 
\newpage

\mynonumsection{Preface}

{\small 
The goal of this paper is to present an algebraic approach to the 
basic results of the theory of linear recurrence relations. 
This approach is based on the ideas from the theory of \emph{representations of one endomorphisms}
(a special case of which may be better known to the reader as the theory of the \emph{Jordan normal form of matrices}\hff).
The notion of the \emph{divided derivatives}, an analogue of the \emph{divided powers}, turned out to be crucial
for proving the results in a natural way and in their natural generality. 
The final form of our methods was influenced by the \emph{the umbral calculus of G.-C. Rota}.

Neither the theory of representation of one endomorphism, nor the theory of divided powers, 
nor the umbral calculus apply directly to our situation.
For each of these theories we need only a modified version of a fragment of it.
This is one of the reasons for presenting all proofs from the scratch.
Both these fragments and our modifications of them are completely elementary and beautiful by themselves.
This is another reason for presenting proofs independent of any advanced sources.
Finally, the theory of the linear recurrence relation is an essentially elementary theory, 
and as such it deserves a self-contained exposition.

The prerequisites for reading this paper are rather modest.
Only the familiarity with the most basic notions of the abstract algebra, such as the notions of
a commutative ring, of a module over a commutative ring, and of endomorphisms and homomorphisms are needed.
No substantial results from the abstract algebra are used. 
A taste for the abstract algebra and a superficial familiarity with it should be sufficient for the reading of this paper.

The standard expositions of the theory of linear recurrence relations present this theory over algebraically closed fields
of characteristic $0$\dfcom or even only over the field of complex numbers. 
In contrast, such restrictions are very unnatural from our point of view. 
The methods of this paper apply equally well to all commutative rings with unit and without zero divisors;
no assumptions about the characteristic are needed.
Of course, a form of the condition of being algebraically closed is needed.
We assume only that all roots of the \emph{characteristic polynomial}\dss
of the linear recurrence relation in question are contained in the ring under consideration
(this can be done using any of the standard approaches to the theory also).
\vspace*{\medskipamount}

\mytitle{Structure of the paper{}.} The main results are stated 
and proved in Section\qss \ref{relations},\qss which depends on all previous ones.
Sections\qss \ref{dder} -- \ref{polynomials}\qss are independent with only one exception:\qss 
Section\qss \ref{operators}\qss depends on both Section\qss \ref{dder}\qss and\qss \ref{seq}.\qss
All references are relegated to the\sss 
\emph{Note bibliographique}\sss at the end.\qss
The reasons are the same as N. Bourbaki's reasons for 
not including any references with the exception of his\sss \emph{Notes historiques}.\qss
\vspace*{\bigskipamount}

}

\renewcommand{\baselinestretch}{1.01}
\selectfont

{\em We denote by\sss $\zz$\sss the ring of integers and by\sss $\nn$\sss the set of non-negative integers.
We denote by\sss $\kk$\sss a fixed entire ring, i.e. a commutative ring with a unit without zero divisors 
and such that its unit is not equal to its zero.\dss
There is a canonical ring homomorphism\sss $\zz\tto\kk$\sss taking\sss $0\dff,\dff 1\elem\zz$ to the zero and the unit of\trs $\kk$\sss respectively, 
making\sss $\kk$\sss into a\sss $\zz$\dnsp-algebra,\sss and every\sss $\kk$\dnsp-module into a\sss $\zz$\dnsp-module.
We identify\sss $0\dff,\dff 1\elem\zz$\sss with their images in\sss $\kk$\dfdot}

\mysection{Divided derivatives of polynomials}{dder}

\myitpar{Polynomials in two variables.} Let $x$ be a variable, 
and let $\kk[x]$ be the $\kk$-algebra of polynomials in $x$ with coefficients in $\kk$\ffdot
Let $y$ be some other variable, and let $\kk[x\fff,y]$ be the $\kk$-algebra of polynomials 
in two variables $x\dff,\dff y$ with coefficients in $\kk$\ffdot
As is well known, $\kk[x\fff,y]$ is canonically isomorphic to the $\kk$-algebra $\kk[x][y]$ of polynomials in $y$ with coefficients in $\kk[x]$\dfdot
We will identify these two algebras.
This allows us to write any polynomial $f(x\fff,y)$ in two variables $x\dff,\dff y$ in the form 
\begin{equation*}
f(x\fff,y)\off =\off \sum\nolimits_{n\eeq 0}^{\infty}\qff g_n(x)\dff y^n\endss.
\end{equation*}
In fact, this sum is obviously finite. 
Equivalently, the polynomials $g_n(x)$ are equal to $0$ for all sufficient big $n\elem\nn$\dfdot
The polynomials $g_n(x)$ are uniquely determined by $f(x\fff,y)$\dfdot

\myitpar{The definition of the divided derivatives.} Let $p(x)\elem\kk[x]$\dfdot
Then $p(x{\pplus}y)\elem\kk[x\fff,y]$\dfcom
and hence $p(x{\pplus}y)$ has the form
\begin{equation}
\label{divided-der}
p(x{\pplus}y)\off =\off \sum\nolimits_{n\eeq 0}^{\infty}\qff (\delta^n\fff p)\fff (x)\dff y^n\endss,
\end{equation}
for some polynomials $(\delta^n\fff p)\fff (x)$
uniquely determined by $p(x)$\dfdot
The sum in\dss (\ref{divided-der})\dss is actually finite. 
Equivalently\hff, $(\delta^n\fff p)\fff (x)\eeq 0$ for all sufficient big $n\elem\nn$\dnsp.\dff\dss
The coefficient $(\delta^n\fff p)\fff (x)$ in front of $y^n$ in the sum in\dss (\ref{divided-der})\dss  
is called \emph{the $n$\dnsp-\emph{th} divided derivative}\dss of the polynomial $p(x)$\dfdot
We will also denote $(\delta^n\fff p)\fff (x)$ 
by $\delta^n\bigl( p\fff (x)\bigr)$ or $\delta^n\fff p (x)$\dfdot

\myitpar{Operators $\delta^n$\dnsp.} Let $n\elem\nn$\dfdot 
By assigning $\delta^n\fff p(x)\elem\kk[x]$ to $p(x)\elem\kk[x]$ we get a map
\[
\delta^n\colon p(x)\qff \longmapsto\qff \delta^n\fff p(x)\endss.
\]
Clearly, $\delta^n$ is a $\kk$\dnsp-linear operator $\kk[x]\to\kk[x]$\dfdot 

After substitution $y\eeq 0$ the equation\dss (\ref{divided-der})\dss reduces to
\[
p(x)\off =\off \delta^0\fff p(x)\endss.
\]
Therefore, $\delta^0\eeq\id\eeq\id_{\kk[x]}$\dnsp.

\mypar{Theorem {(}{\fff}Leibniz formula\fff{)}.}{leibniz-formula} 
\emph{Let $f(x)\dff,\dff g(x)\elem\kk[x]$\dfcom and let $n\elem\nn$\dfdot
Then}
\begin{equation*}
\delta^n \left( f(x) g(x)\right)\off 
=\off \sum\nolimits_{\fff i{\pplus}j\eeq n}\qff \delta^{\fff i}\fff f(x)\dff \delta^{\fff j} g(x)\endss.
\end{equation*}

\newpage
\proof By applying\dss (\ref{divided-der})\dss to $p(x)\eeq f(x)$ and to $p(x)\eeq g(x)$\dfcom
we get
\begin{equation*}
f(x{\pplus}y)\off 
=\off \sum\nolimits_{\fff i\eeq 0}^{\fff \infty}\qff \delta^{\fff i}\fff f\fff (x)\dff y^{\fff i}\endss,
\end{equation*}
\begin{equation*}
g(x{\pplus}y)\off 
=\off \sum\nolimits_{\fff j\eeq 0}^{\fff \infty}\qff \delta^{\fff j}\fff g\fff (x)\dff y^{\fff j}\endss.
\end{equation*}

\vspace*{-\bigskipamount}
By multiplying these two identities, we get 
\begin{equation*}
f(x{\pplus}y)\dff g(x{\pplus}y)\off 
=\off \sum\nolimits_{\fff i,j\eeq 0}^{\fff \infty}\qff 
\delta^{\fff i}\fff f\fff (x)\dff y^{\fff i}\dff\dff \delta^{\fff j}\fff g\fff (x)\dff y^{\fff i}\endss,
\end{equation*}
and hence
\begin{align*}
f(x{\pplus}y)\dff g(x{\pplus}y)\off &=\off \sum\nolimits_{\fff i,j\eeq 0}^{\fff \infty}\qff 
\delta^{\fff i}\fff f\fff (x)\dff\dff \delta^{\fff j}\fff g\fff (x)\dff y^{i{\pplus}j} \\
                                &=\off \sum\nolimits_{\fff n\eeq 0}^{\fff \infty}\qff \left(\hspace*{0.05em} \sum\nolimits_{i{\pplus}j\eeq n}\qff
                                \delta^{\fff i}\fff f\fff (x)\dff\dff \delta^{\fff j}\fff g\fff (x)\right) y^{n}\endss.
\end{align*}
The theorem follows.  \eproof

\mypar{Corollary {(}{\fff}Leibniz formula for $\delta^1$\dnsp{)}\hff.}{leibniz-formula-one} 
\emph{Let $f(x)\dff,\dff g(x)\elem\kk[x]$\dfdot
Then}
\begin{equation*}
\delta^1\left(f(x) g(x)\right)\off =\off \delta^1\fff f(x)\dff g(x)\qff +\qff f(x)\dff \delta^1 g(x)\endss.
\end{equation*}
\emph{In other terms, $\delta^1$ is a \textbf{derivation} of the ring $\kk[x]$\dnsp.}  \eproof

\mypar{Lemma.}{der-of-xn01} (i) \emph{$\delta^0\fff (1)\qff =\qff 1$ 
and $\delta^n\fff (1)\qff =\qff 0$\sss for all{\fff}\dss $n\geqs 1$\dfdot} 
\vspace*{\medskipamount} \\
\hspace*{5.54em} (ii) \emph{$\delta^0\fff x\qff =\qff x$\dfcom 
$\delta^1\fff x\qff =\qff 1$\dfcom and $\delta^n\fff x\qff =\qff 0$ for all\dss $n\geq 2$\dfdot}

\proof As we noted above, $\delta^0\eeq\id$\dfdot 
In particular, $\delta^0\fff (1)\eeq 1$\dfdot
For $p(x)\eeq 1$\dfcom the formula (\ref{divided-der}) takes the form
\begin{equation}
\label{delta1}
1\off =\off \delta^0\fff (1)\dff y^0\qff +\qff  \sum\nolimits_{n\eeq 1}^{\infty}\qff \delta^n\fff (1)\dff y^n\endss.
\end{equation}
Since $\delta^0\dff (1)\fff y^0 = 1\cdot y^0 = 1$\dfcom the formula (\ref{delta1}) implies that
\begin{equation*}
0\off =\off  \sum\nolimits_{n\eeq 1}^{\infty}\qff \delta^n\fff (1)\dff y^n\endss.
\end{equation*}
It follows that  $\delta^n\fff 1\eeq 0$ for $n\geqs 1$\dfdot
This proves the part\sss (i)\sss of the lemma.
For $p(x)\eeq x$\dfcom the formula\sss (\ref{divided-der})\sss takes the form
\[
x + y\off =\off \delta^0\fff (x)\dff y^0\qff +\qff 
\delta^1\fff (x)\dff y^1\qff +\qff \sum\nolimits_{n\eeq 2}^{\infty}\qff \delta^n\fff (x)\dff y^n\endss.
\]
It follows that $\delta^0\fff (x)\eeq x$\dfcom  $\delta^1\fff (x)\eeq 1$\dfcom and $\delta^n\fff (x)\eeq 0$ for all\dss $n\geqs 2$\dfdot
This proves the part\sss (ii)\sss of the lemma.  \eproof \newpage

\mypar{Lemma.}{delta-one-of-xn} $\delta^1\fff (x^n)\off =\off n\fff x^{n{\minus}1}$ 
{\dff}\emph{for all}\trs $n\elem\nn$\dfcom $n\geq 1$\dfdot

\proof The case $n\eeq 1$ was proved in Lemma \ref{der-of-xn01}.  
Suppose that $n\elem\nn$\dnsp,\dff\dss $n\geqs 1$\dfcom 
and we already know that $\delta^1 (x^n)\eeq n\fff x^{n{\minus}1}$\dfdot
By Corollary \ref{leibniz-formula-one},\vspace{\medskipamount}
\begin{align*}
\delta^1\fff (x^{n{\plus}1})\off &=\off \delta^1\fff (x\cdot x^n)\off 
=\off \delta^1\fff (x)\fff x^n\qff +\qff x\dff \delta^1\fff (x^n) \\ 
                             &=\off 1\cdot x^n\qff +\qff x\dff (n\fff x^{n{\minus}1})\off 
                             =\off x^n\qff +\qff n\fff x^n\off 
                             =\off (n{\plus}1)\dff x^n\endss.
\end{align*}

\vspace{-\medskipamount}
An application of induction completes the proof.  \eproof

\myitpar{Remark.} By Lemma \ref{delta-one-of-xn}, the operator $\delta^1\colon\kk[x]\tto\kk[x]$
agrees on the powers $x^n\elem\kk[x]$ with the operator $d\colon f(x)\mapsto f'(x)$ of taking the usual formal derivative.
Since both these operators are  $k$\dnsp-linear, $\delta^1\eeq d$\dfdot
But if $i\elem\nn$\dfcom $i\geqs 2$\dfcom
then the operator $\delta^{\fff i}$ is \emph{not equal}\dss to the operator of taking the $i$\dnsp-th derivative.
This immediately follows either from Lemma \ref{dder-xn} or from Theorem \ref{composition} below.

\myitpar{Binomial coefficients.}  Let $n\elem\nn$\dfdot
For $i\elem\nn$\dfcom $i\leqs n$\dfcom
we\dss \emph{define}\dss the\dss \emph{binomial coefficients}\dss $\bico{i{\mminus}n}{n}\elem\nn$\dnsp,\dss 
by the\sss \emph{binomial formula}
\begin{equation}
\label{binomial-formula}
(x + y)^{i}\off 
=\off \sum\nolimits_{\fff n\eeq 0}^{\fff i}{\dff}\qff 
\bico{i{\mminus}n}{n}\dff\dff x^{\fff i{\mminus}n\dff}\fff y^{\fff n}\endss.
\end{equation}
Given arbitrary numbers $a\dff,\dff b\elem\nn$\dfcom 
we define $\bico{a}{b}$ as $\bico{n{\mminus}b}{b}$\dfcom
where $n\eeq a{\pplus}b$\dnsp.\dss
Given arbitrary integers $a\dff,\dff b\elem\zz$\dfcom 
we set $\bico{a}{b}\eeq 0$ if at least one of the numbers $a$\dfcom $b$ is not in $\nn$\dfdot

We prefer the notation $\bico{a}{b}$ to the classical one by the typographical reason,
and because the new notation helps to bring to the light the fact that we do not use any properties of
$\bico{a}{b}$ except the above definition.

\mypar{Lemma.}{dder-xn} $\delta^{\fff n}\dff (x^{\fff i})\off 
=\off \bico{i{\mminus}n}{n}\dff\dff x^{\fff i{\mminus}n}$\dfdot

\proof It is sufficient to compare\dss (\ref{binomial-formula})\dss 
with the definition\dss (\ref{divided-der})\dss of divided derivatives. \eproof

\myitpar{The left shift operator.} The\qss \emph{left shift operator} $\lambda\colon\kk[x]\tto\kk[x]$ is just the operator of multiplication by the polynomial $x$\nsp:
\[
\lambda\dff p(x)\off =\off \lambda (p) (x)\off =\off x\dff p(x)\endss.
\]\newpage
The reasons for calling $\lambda$ the\qss \emph{left shift operator}\qss will be clear later.\qss
The main property of the left shift operator and the divided derivatives is the following commutation relation.

\mypar{Theorem.}{commutator} \emph{Let us set $\delta^{{\minus}1}\eeq 0$\dfdot
Then}
\begin{equation}
\label{delta-lambda}
\delta^n\circ\lambda\qff -\qff \lambda\circ\delta^n\off =\off \delta^{n{\minus}1}\endss.
\end{equation}
\emph{for all}\sss $n\elem\nn$\dfdot

\proof Let $p(x)\elem\kk[x]$\dfcom and let $y$ be a variable different from $x$\dfdot
By the Leibniz formula from Theorem \ref{leibniz-formula},
\[
\delta^n\fff \bigl(x\dff p(x)\bigr)\off 
=\off \sum\nolimits_{i{\pplus}j\eeq n}\qff \delta^{\fff i}\fff (x)\dff\dff \delta^{\fff j}\fff\bigl( p(x)\bigr)\endss.
\]
But by Lemma \ref{der-of-xn01}, $\delta^0\fff x\eeq x$\dfcom $\delta^1\fff x\eeq 1$\dfcom 
and $\delta^{\fff i}\fff x\eeq 0$ for $i\geqs 2$\dfdot
Therefore
\[
\delta^n\fff \bigl(x\dff p(x)\bigr)\off 
=\off x\dff\delta^n\fff\bigl( p(x)\bigr)\qff 
+\qff \delta^{n{\minus}1}\fff\bigl( p(x)\bigr)\endss,\hspace*{0.5em}\mbox{ or,\dss what is the same, }
\]
\[
\delta^n\fff (x\dff p(x))\qff -\qff x\dff\delta^n\fff\bigl( p(x)\bigr)\off 
=\off \delta^{n{\minus}1}\fff\bigl( p(x)\bigr)\endss.
\]

\vspace*{-\medskipamount}
Rewriting the last identity in terms of $\lambda$\dfcom we get
\[
\delta^n\fff \bigl(\fff\lambda (p) (x)\bigr)\qff -\qff \lambda\fff\bigl(\delta^n\bigl( p(x)\bigr)\bigr)\off
=\off \delta^{n{\minus}1}\fff\bigl( p(x)\bigr),\hspace*{0.4em}\mbox{ i.e. }
\]
\[
\delta^n\circ\lambda \bigl(p(x)\bigr)\qff -\qff \lambda\circ\delta^n\fff\bigl(\fff p(x)\bigr)\off 
= \off\delta^{n{\minus}1}\bigl(p(x)\bigr)\endss.
\]

\vspace*{-\medskipamount}
Since $p(x)\elem\kk[x]$ was arbitrary, this proves the theorem.  \eproof

\myitpar{The following theorem will be not used in the rest of the paper.}

\mypar{Theorem (Composition of divided derivatives).}{composition} \emph{Let $n\dff,\dff m\elem\nn$\dfdot
Then}
\[
\dff\delta^n\circ\delta^m\off =\off \bico{m}{n}\dff\dff \delta^{n+m} \endss.
\]

\vspace*{-\bigskipamount}
\proof Let $p(x)\elem\kk[x]$\dfdot
Let $u\dff,\dff z$ be two new variables different from both $x$ and $y$\dfdot
If we apply\sss (\ref{divided-der})\sss to $u\dff,\dff z$ in the role of $x\dff,\dff y$ respectively (and use $m$ instead of $n$\dnsp), we get
\begin{equation*}
p(u{\pplus}z)\off =\off \sum\nolimits_{m\eeq 0}^{\infty}\qff (\delta^m\fff p)\fff (u)\dff z^m\endss.
\end{equation*}
Let us set $u=x{\pplus}y$ and apply\sss (\ref{divided-der})\sss 
to each polynomial $(\delta^m\fff p)\fff (x{\pplus}y)$\nsp:\newpage
\begin{align}
\label{xy-z}
p((x{\pplus}y){\pplus}z)\off &=\off \sum\nolimits_{m\eeq 0}^{\infty}{\dff}\qff (\delta^m\fff p)\fff (x{\pplus}y)\dff\dff z^m\\
                             &=\off \sum\nolimits_{m\eeq 0}^{\infty}{\dff}\qff 
                             \left(\sum\nolimits_{n\eeq 0}^{\infty}\qff \delta^n\fff (\delta^m\fff p)\fff (x)\dff\dff y^n\right) z^m \notag \\
                                     &=\off \sum\nolimits_{m\dff,\dff n\eeq 0}^{\infty}\qff \delta^n\circ\delta^m\fff (p)\fff (x)\dff\dff\dff y^n\dff z^m \notag                  \end{align}
Alternatively, we can apply\sss (\ref{divided-der})\sss to $y{\pplus}z$ in the role of $y$
and then apply\sss (\ref{binomial-formula}):
\begin{align}
\label{x-yz}
p(x{\pplus}(y{\pplus}z))\off &=\off \sum\nolimits_{k\eeq 0}^{\infty}{\dff}{\dff} \delta^k\dff p\fff (x)\dff\dff\dff (y{\pplus}z)^k \\
                         &=\off \sum\nolimits_{k\eeq 0}^{\infty}{\dff}{\dff} \delta^k\dff p\fff (x)\dff\left(\sum\nolimits_{m+n=k}{\qff}\bico{m}{n}\dff  y^m\fff z^n\right) \notag \\ 
\phantom{p(x{\pplus}(y{\pplus}z))}\off &=\off \sum\nolimits_{m\dff,\dff n\eeq 0}^{\infty}{\dff}\qff \delta^{m+n}\dff p\fff (x)\dff\dff\dff \bico{m}{n}\dff  y^m\fff z^n \notag \\
\phantom{p(x{\pplus}(y{\pplus}z))}\off &=\off \sum\nolimits_{m\dff,\dff n\eeq 0}^{\infty}{\dff}\qff  \bico{m}{n}\dff \delta^{m+n}\dff p\fff (x)\dff\dff\dff  y^m\fff z^n\endss. \notag 
\end{align} 

\vspace*{-\medskipamount}
By the associativity of the addition, $(x{\pplus}y){\pplus}z\eeq x{\pplus}(y{\pplus}z)$ and hence
\[
p((x{\pplus}y){\pplus}z)\off =\off p(x{\pplus}(y{\pplus}z))\endss.
\]
By combining this equality with\sss (\ref{xy-z})\sss and\sss (\ref{x-yz})\sss
we conclude that
\[
\delta^n\circ\delta^m\fff (p(x))\off =\off \bico{m}{n}\dff (\delta^{m+n}\dff p(x))
\]
for all $p(x)\elem\kk[z]$ and $n\dff,\dff m\elem\nn$\dfdot
The theorem follows.  \eproof

\mysection{Sequences and duality}{seq}

\myitpar{Sequences.} A \emph{sequence} of elements of a set $X$ is defined as a map $\nn\tto X$\dfdot
For a sequence $s$ we usually denote the value $s(i)\dff,\dff i\elem\nn$ by $s_i$ and often call it the $i$\dnsp-th\sss \emph{term of} $s$\dfdot
The set of all sequences of elements of $X$ will be denoted by $\seqq_X$\dfdot
We are, first of all, interested in the case when $X$ is a $\kk$\dnsp-module, 
and especially in the case when $X$ is equal to $\kk$ considered as a $\kk$\dnsp-module.
When it is clear from the context to what set $X$ the terms of the considered sequences belong, 
we call the sequences of elements of $X$ simply \emph{sequences}.

Let $M$ be a $\kk$\dnsp-module. Then the set $\seqq_M$ has a canonical structure of a $\kk$\dnsp-module.
The $\kk$\dnsp-module operations on $\seqq_M$ are the term-wise addition of sequences and the term-wise multiplication 
of sequences by elements of $\kk$\ffcom defined in the following obvious way{}.
The\sss \emph{term-wise sum}\sss $r{\pplus}s$ of sequences $r\dff,\dff s\elem\seqq_M$ is defined by $(r{\pplus}s)_i\eeq r_i{\pplus}s_i$\ffcom
and the \emph{term-wise product}\sss $c\fff s$ of $c\elem\kk$ and $s\elem\seqq_M$ is defined by $(c\fff s)_i\eeq c\hff s_i$\ffdot

\myitpar{Modules of homomorphisms.} Let $M'\dff,\dff M''$ be $\kk$\dnsp-modules. 
Then the set $\hhom (M'\dff,\dff M'')$ of $\kk$\dnsp-homomorphisms $M'\tto M''$ has a canonical structure of a $\kk$\dnsp-module.
The addition is defined as the addition of $\kk$\dnsp-homomorphisms, and the product $a\hff F$ of an element $a\elem\kk$ and 
a $\kk$\dnsp-homomorphisms $F\colon M'\tto M''$ is defined by $(a\hff F)(m)\eeq a\hff F(m)$\dfcom where $m\elem M'$\dfdot
Note that the (obvious) verification of the fact that $a\hff F$ is a $\kk$\dnsp-homomorphism uses the commutativity of $\kk$\ffdot

We are mostly interested in the case of $M'\eeq\kk[x]$\dfcom 
and especially in the case of $M'\eeq\kk[x]$ and  $M''\eeq\kk$\dfcom where $\kk[x]$ is considered as a $\kk$\dnsp-module 
by forgetting about the multiplication of elements of $\kk[x]$\dfcom 
and the ring $\kk$ is considered as a module over itself.

\vspace*{\bigskipamount}
\textbf{\textit{For the rest of this section\dss $M$\dss denotes a fixed\dss $\kk$\dnsp-module.}}

\myitpar{A pairing between $\seqq_M$ and $\kk[x]$\dnsp.} Consider a sequence $s\elem\seqq_M$ and a polynomial
\[
p(x)\off =\off \sum\nolimits_{i\eeq 0}^{\infty}\qff c_i\dff x^i\in\kk[x]\endss.
\]
Of course, the sum here is actually  finite, i.e. $c_i\eeq 0$ for all sufficiently large $i$\dfdot
Let
\[
\aang{s}{p(x)}\off =\off \sum\nolimits_{i\eeq 0}^{\infty}\qff c_i\dff s_i\in M\endss.
\]
Since $c_i\eeq 0$ for all sufficiently large $i$\dfcom the sum in the right hand side of this formula is well defined.
The map 
\[
\aang{\bullet\dff}{\bullet}\colon\seqq_M\times\kk[x]\qff \longrightarrow\qff M
\]
defined by 
\[
\aang{\bullet\dff}{\bullet}\colon (s\dff,\dff p(x))\qff \longmapsto\qff \aang{s}{p(x)}
\]
is our \emph{pairing between $\kk[x]$ and} $\seqq_M$\dfdot
Obviously, it is a $\kk$\dnsp-bilinear map (and hence indeed deserves to be called a \emph{pairing}{\fff}).

The pairing $\aang{\bullet\dff}{\bullet}$ defines a $\kk$\dnsp-linear map
\[
\dual_M\colon\seqq_M\qff \longrightarrow\qff \hhom \bigl(\hff \kk[x]\dff,\dff M\bigr)
\]
by the usual rule $\dual_M (s)\colon p(x)\mapsto \aang{s}{p(x)}$\dfdot 

Note that, obviously, $\aang{s}{x^i} = s_{\fff i}$
for every $s\elem\seqq_M$\dfcom $n\elem\nn$\dfdot
Therefore
\begin{equation}
\label{powers-dual}
s_{\fff i}\off =\off \dual_M (s) \bigl( x^i \bigr)
\end{equation}
for every $i\elem\nn$ and every $s\elem\seqq_M$\dfdot

\mypar{Theorem ({\fff}Duality{\fff}).}{duality} \emph{The pairing $\aang{\bullet\dff}{\bullet}$ is non-degenerate in the sense that
the map}
\[
\dual_M\colon\seqq_M\qff \longrightarrow\qff \hhom \bigl(\hff \kk[x]\dff,\dff M \bigr)
\]
\emph{is an isomorphism of $\kk$\dnsp-modules.}

\proof Note that $\kk$\dnsp-homomorphism $F\colon\kk[x]\tto M$ is 
determined by its values $F(x^i)$ on the monomials $x^i\dff,\dff i\elem\nn$
({\fff}because every polynomial $p(x)\elem\kk[x]$ is a finite sum of 
powers $x^i\dff,\dff i\elem\nn$ with coefficients in $\kk$).
By (\ref{powers-dual}) the terms $s_i$ of a sequence $s\elem\seqq_M$ are equal to the values $\dual_M (s) (x^i)$\dfdot
It follows that all terms of $s$\dfcom and hence the sequence $s$ are determined by the homomorphism  $\dual_M (s)$\dfdot
Therefore, the map $\dual_M$ is injective.

In order to prove that $\dual_M$ is surjective, 	
let us consider an arbitrary $\kk$\dnsp-homomorphism $F\colon\kk[x]\tto M$\dnsp.\dss
Let $s\elem\seqq_M$ be the sequence with $s_i\eeq F(x^i)$\dfdot
Then $\kk$\dnsp-homomorphisms $F$ and $\dual_M (s)$ take the same values at all powers $x^i$\dfdot
It follows that $F\eeq\fff\dual_M (s)$ (cf. the previous paragraph).
Therefore, the map $\dual_M$ is surjective.\qss
The theorem follows.  \eproof

\myitpar{Dual endomorphisms.} Each $\kk$\dnsp-endomorphism $E\colon\kk[x]\tto\kk[x]$ 
defines\dss \emph{its dual endomorphism} 
\[
E^*\colon\hhom \bigl(\fff\kk[x]\dff,\dff M\bigr)\qff \longrightarrow\qff \hhom \bigl(\fff\kk[x]\dff,\dff M\bigr)
\]
by the usual formula $E^*(h)\eeq h\circ E$\dfcom
for all $\kk$\dnsp-homomorphisms $h\colon\kk[x]\tto M$\dfdot 
Obviously, if $E\dff,\dff F$ are two $\kk$\dnsp-endomorphisms $\kk[x]\tto\kk[x]$\dfcom
then $(E\circ F)^* = F^*\circ E^*$\dfdot

\myitpar{Adjoint endomorphisms.} Let $M$ be a $\kk$\dnsp-module.
Since $\dual_M$ is an isomorphism by Theorem \ref{duality}, 
we can use $\dual_M$ to turn the dual map 
\[
E^*\colon\hhom (\kk[x]\dff,\dff M)\qff \longrightarrow\qff \hhom (\kk[x]\dff,\dff M)
\]
of an endomorphism $E\colon\kk[x]\tto\kk[x]$ into a map $\seqq_M\tto\seqq_M$\dfdot
Namely, let 
\[
E^{\perp}\off =\off (\dual_M)^{{\minus}1}\circ E^*\circ\dual_M\colon\seqq_M\qff \longrightarrow\qff \seqq_M\endss.
\]
Then $\dual_M\circ E^{\perp} = E^{*}\circ\dual_M$\dfdot
We will call $E^\perp$ the\qss \emph{adjoint endomorphism of}\qss $E$\dfdot
\newpage

For every pair $E\dff,\dff F\colon\kk[x]\tto\kk[x]$ of $\kk$\dnsp-endomorphisms $(E\circ F)^\perp = F^\perp\circ E^\perp$\dfdot
This immediately follows from the corresponding property $(E\circ F)^* = F^*\circ E^*$ of dual endomorphisms.

\mypar{Lemma.}{adjoint-seq} \emph{Let $M$ be a $\kk$\dnsp-module,
and let $E\colon\kk[x]\tto\kk[x]$ be a $\kk$\dnsp-endomorphism.
The adjoint map $\seqq_M\tto\seqq_M$ is the unique map $E^\perp$ such that}
\begin{equation}
\label{perp-adj}
\aang{p}{E^\perp(s)}\off =\off \aang{E(p)}{s}
\end{equation}
\emph{for all $p\eeq p(x)\elem\kk[x]$\dfcom $s\elem\seqq_M$\dfdot}

\proof Let $p\eeq p(x)\elem\kk[x]$ and let $s\elem\seqq_M$\dfdot
By the definition of $\dual_M$ we have: 
\begin{align*}
\aang{p}{E^\perp(s)}\off &=\off \dual_M \bigl(\fff E^\perp (s)\fff\bigr) \bigl( p\bigr)\off 
=\off \bigl(\fff\dual_M\circ E^\perp (s)\fff\bigr)\bigl( p\bigr)\endss; \\
\aang{E(p)}{s}\off &=\off \dual_M (s) \bigl(\fff E(p)\bigr)\off =\off E^* \bigl(\fff\dual_M (s) \fff\bigr)\bigl( p\bigr)\off 
=\off \bigl( E^*\circ\dual_M (s)\bigr)\bigl( p\bigr)\endss. 
\end{align*}
Therefore, (\ref{perp-adj}) is equivalent to
$\displaystyle
\bigl(\fff\dual_M\circ E^\perp (s)\fff\bigr)\bigl( p\bigr) 
=\off \bigl( E^*\circ\dual_M (s)\bigr)\bigl( p\bigr)$\dfdot

It follows that (\ref{perp-adj}) holds for all $p\eeq p(x)\elem\kk[x]$\dfcom $s\elem\seqq_M$
if and only if $\dual_M\circ E^{\perp} =\linebreak E^{*}\circ\dual_M$\dfdot
The lemma follows.  \eproof

\mysection{Adjoints of the left shift and of divided derivatives}{operators}

\vspace*{\bigskipamount}
\textbf{\textit{As in the previous section,\dss $M$\dss denotes a fixed\sss $\kk$\dnsp-module.}}

\myitpar{The adjoint of the left shift operator.} Let $L\eeq\lambda^\perp$\dfcom where 
$\lambda$ 
is the left shift operator from Section \ref{dder}.\dss
In view of the following lemma call $L$ also the left shift operator.

\mypar{Lemma.}{left-shift} \emph{For every sequence $s\elem\seqq_M$ the terms of the sequence $L\fff (s)$
are $\bigl(L\fff (s)\bigr)_{i}\eeq s_{\fff i{\plus}1}$\dfdot}

\proof Recall that $s_{\fff i}\qff =\qff \aang{s}{x^i}$ for any sequence $s\elem\seqq_M$\dfdot 
Therefore by Lemma \ref{adjoint-seq}
\begin{equation*}
\bigl(L\fff (s)\bigr)_i\off 
=\off \aang{L\fff (s)}{x^i}\off 
=\off \aang{\lambda^\perp(s)}{x^i}\off 
                          =\off \aang{s}{\lambda\fff (x^i)}\off 
                          =\off \aang{s}{x^{i{\plus}1}}\off 
                          =\off s_{\fff i{\plus}1}\endss.
\end{equation*}
The lemma follows.  \eproof \newpage

\mypar{Corollary.}{iterated-shift} \emph{For every $n\elem\nn$ and every sequence $s\elem\seqq_M$ 
the terms of the sequence $L^n\fff (s)$
are $\bigl(L^n\fff (s)\bigr)_{i}\qff =\qff s_{\fff i{\plus}n}$\dfdot}  

\proof For $n\eeq 0$ the corollary is trivial, because $L^0\eeq\id$\dfdot
For $n\geqs 1$ the corollary follows from Lemma \ref{left-shift}, if we use an induction by $n$\dfdot  \eproof

\myitpar{The adjoints of the divided derivatives.} Let $\deseq{n}\eeq (\fff\delta^n\fff)^\perp$\dnsp,\dss where $n\elem\nn$ and
$\delta^n$ 
is the $n$\dnsp-th divided derivative operator from Section \ref{dder}.
Recall (see Section \ref{dder}) that $\delta^0\colon\kk[x]\tto\kk[x]$ is the identity of $\kk[x]$\dfdot
Therefore $D^0\colon\seqq_M\tto\seqq_M$ is also the identity of $\seqq_M$\dfdot

Recall that in Theorem \ref{commutator} we also introduced operator $\delta^{{\minus}1}$\dfdot
Let $\deseq{{\minus}1}\eeq (\delta^{{\minus}1})^\perp$\dfdot
Since $\delta^{{\minus}1}\eeq 0$ by the definition, we have $\deseq{{\minus}1}\eeq 0$\dfdot

The following commutation relations are the most important for our purposes 
properties of the adjoint operators $L\eeq\lambda^\perp$ and $\deseq{n}\eeq(\delta^n)^\perp$\dfdot

\mypar{Theorem.}{seq-commutator} \emph{For every $\alpha\elem\kk$ and every $n\elem\nn$}
\begin{equation*}
L\circ\deseq{n}\qff -\qff \deseq{n}\circ L\off 
=\off \deseq{n{\minus}1}\hspace*{0.35em}\mbox{ \emph{and} }\hspace*{0.75em}
\end{equation*}
\begin{equation*}
(L{\mminus}\alpha)\circ\deseq{n}\qff -\qff \deseq{n}\circ (L{\mminus}\alpha)\off 
=\off \deseq{n{\minus}1}\nsp, 
\end{equation*}

\vspace*{-\medskipamount}
\emph{where we interpret $\alpha$ as the operator $\seqq_M\qff \longrightarrow\qff \seqq_M$ 
of multiplication by $\alpha\elem\kk$\ffdot}

\proof By taking the adjoint identity of the identity\dss 
(\ref{delta-lambda})\dss from Theorem \ref{commutator}\fff, we get
\[
\bigl(\delta^n\circ\lambda\bigr)^\perp\qff -\qff \bigl(\lambda\circ\delta^n\bigr)^\perp\off 
= \off\bigl(\delta^{n{\minus}1}\bigr)^\perp,
\]
and hence 
\[
\lambda^\perp\circ\bigl(\delta^n\bigr)^\perp\qff -\qff \bigl(\delta^n\bigr)^\perp\circ\lambda^\perp\off 
=\off \bigl(\delta^{n{\minus}1}\bigr)^\perp.
\]
In view of the definitions of $L$ and $D^n$\dfcom this implies the first identity of the theorem. 

Since $\deseq{n}$ is a $\kk$\dnsp-linear operator, we have $\alpha\circ\deseq{n}\eeq\deseq{n}\circ\alpha$,
where $\alpha$ is interpreted as the multiplication operator.
Clearly, the first identity of the theorem together with $\alpha\circ\deseq{n}\eeq\deseq{n}\circ\alpha$ implies the second one.  \eproof

\myitpar{The sequences\sss $s(\alpha)$\sss and\sss $s(\alpha\dff,\dff n)$\dnsp.} Let $\alpha\elem\kk$ and $n\elem\nn$\dfdot
Let us define sequences $s(\alpha)$ and $s(\alpha\dff,\dff n)$ by
\[
s(\alpha)_{\fff i}\off =\off \alpha^{\fff i}
\hspace*{0.8em}\mbox{ and }\hspace*{0.8em}
s(\alpha\dff,\dff n)\off =\off \deseq{n}\bigl(s(\alpha)\bigr)\endss.
\]
Obviously, $s(\alpha\dff,\dff 0)\eeq s(\alpha)$\dfdot
Note that $s(\alpha)\nneq 0$ even if $\alpha\eeq 0$\dnsp,\dss 
because $s(\alpha)_0\eeq\alpha^0\eeq 1$ by the definition for all $\alpha\elem\kk$\dfdot

For explicit formulas for sequences $s(\alpha\dff,\dff n)$ with $n\geqs 1$ 
the reader is referred to Theorem \ref{formulas-delta-s}\sss below.
No such formulas are used in this paper.

\mypar{Lemma.}{eigen-powers} \emph{Let\sss $\alpha\elem\kk$\dss and\sss $s\elem\seqq_M$\dfdot
Then $\bigl(L{\mminus}\alpha\bigr) (s) = 0$ if and only if $s$ has the form $s\eeq \beta\fff s(\alpha)$\dfcom
where $\beta\elem\kk$\dfdot}

\proof The condition $\bigl(L{\mminus}\alpha\bigr) (s)$ is equivalent to\dss
$L(s)\eeq\alpha\fff s$\dfdot
The latter condition holds if and only if 
$\bigl(L\fff (s)\bigr)_{n} = \alpha\dff s_{\fff n}$ for all $n\elem\nn$\dfdot 
By Lemma \ref{left-shift} $\bigl(L\fff (s)\bigr)_{n}\eeq s_{\fff n{\plus}1}$\dnsp.\dss
Therefore, $\bigl(L{\mminus}\alpha\bigr) (s)\eeq 0$ if and only if
$s_{\fff n{\plus}1} = \alpha\dff s_{\fff n}$ for all $n\elem\nn$\dfdot
An application of the induction completes the proof.  \eproof

\mypar{Lemma.}{left-der} \emph{Suppose that\sss $a\elem\nn$\dfcom $\alpha\elem\kk$\dfcom and\sss $s\elem\seqq_M$\dfdot
If\trs $n\geqs a$\dfcom then} 
\[
\bigl(L{\mminus}\alpha\bigr)^a\bigl(s(\alpha\dff,\dff n)\bigr)\off 
=\off s(\alpha\dff,\dff n{\minus}a)\endss.
\]

\vspace*{-\bigskipamount}
\proof The lemma is trivial if $a\eeq 0$\dfdot 
Let us prove the lemma for $a\eeq 1$\dfdot
In view of the definition of sequences $s(\alpha\dff,\dff n)$\dfcom we need to prove that 
\[
\bigl(L{\mminus}\alpha\bigr)\bigl(\deseq{n}\bigl( s(\alpha)\bigr)\bigr)\off 
=\off \deseq{n{\minus}1}\bigl( s(\alpha)\bigr)\endss.
\]
By applying the second identity of Theorem \ref{seq-commutator}\dss to $s(\alpha)$\dfcom we get
\[
(L{\mminus}\alpha)\circ\deseq{n}\bigl( s(\alpha)\bigr)\qff -\qff \deseq{n}\circ (L{\mminus}\alpha) \bigl( s(\alpha)\bigr)\off 
=\off \deseq{n{\minus}1} \bigl( s(\alpha)\bigr)\endss,
\]
which is equivalent to
\[
(L{\mminus}\alpha)\bigl(\deseq{n}\bigl( s(\alpha)\bigr)\bigr)\qff 
-\qff \deseq{n}\bigl( (L{\mminus}\alpha) \bigl( s(\alpha)\bigr)\bigr)\off 
=\off \deseq{n{\minus}1} \bigl( s(\alpha)\bigr)\endss.
\]
Since $\bigl(L{\mminus}\alpha\bigr) \bigl(s(\alpha)\bigr) = 0$ by Lemma \ref{eigen-powers}, we see that
\[
(L{\mminus}\alpha)\bigl(\deseq{n}\bigl( s(\alpha)\bigr)\bigr)\off 
=\off \deseq{n{\minus}1} \bigl( s(\alpha)\bigr)\endss,
\]
i.e. $\bigl(L{\mminus}\alpha\bigr)\bigl(s(\alpha\dff,\dff n)\bigr)\qff =\qff s(\alpha\dff,\dff n{\minus}1)$\dfdot  
This proves the lemma for $a\eeq 1$\dfdot 
The general case follows from this one by induction. \eproof

\myitpar{The following theorem is not used in the rest of the paper.}\newpage

\mypar{Theorem.}{formulas-delta-s} \emph{Let\dss $n\elem\nn$\dfdot
For every sequence\sss $s\elem\seqq_M$\sss the terms of the sequence\sss $\deseq{n} (s)$\sss
are} 
\[
\bigl(\deseq{n} (s)\bigr)_{i}\off =\off \bico{i{\mminus}n}{n}\dff s_{\fff i{\mminus}n}\endss.
\]
\emph{In addition, for every\sss $\alpha\elem\kk$\sss the terms of the sequence\sss $s(\alpha\dff,\dff n)$\sss are} 
\[
\bigl(s(\alpha\dff,\dff n)\bigr)_{\fff i}\off =\off \bico{i{\mminus}n}{n}\dff \alpha^{\fff i{\minus}n}\endss.
\]

\proof Recall that $s_{\fff i}\qff =\qff \aang{s}{x^{\fff i}}$ for any sequence $s\elem\seqq_M$\dfdot 
Together with Lemma \ref{adjoint-seq} this fact implies that 
\begin{align*}
\bigl(\deseq{n} (s)\bigr)_{\fff i}\off &=\off \aang{\deseq{n} (s)}{x^{\fff i}} \\ 
                                   &=\off \aang{(\delta^n)^\perp(s)}{x^{\fff i}} \\ 
                                   &=\off \aang{s}{\delta^{\fff n} (x^{\fff i})} \endss. 
\end{align*}
Since $\delta^{\fff n}\dff (x^{\fff i})\qff 
=\qff \bico{i{\mminus}n}{n}\dff\dff x^{\fff i{\mminus}n}$ by Lemma \ref{dder-xn}, we have
\begin{align*}
\aang{s}{\deseq{n} (x^i)}\off &=\off \aang{s}{\bico{i{\mminus}n}{n}\dff\dff x^{i{\mminus}n}}  \\
                          &=\off \bico{i{\mminus}n}{n}\dff\aang{s}{\deseq{n} (x^{\fff i{\mminus}n})} \\
                          &=\off \bico{i{\mminus}n}{n}\dff s_{\fff i{\mminus}n}\endss.
\end{align*}
The first part of the theorem follows.
Let us apply the first part to $s\eeq s(\alpha)$\dfdot
We get
\begin{align*}
s(\alpha\dff,\dff n)\off &=\off \bigl(\deseq{n} \bigl(s(\alpha)\bigr)\bigr)_{i} \\
                     &=\off \bico{i{\mminus}n}{n}\dff s(\alpha)_{\fff i{\mminus}n}\off 
                     =\off \bico{i{\mminus}n}{n}\dff \alpha^{\fff i{\mminus}n}\endss.
\end{align*}
This proves the second part of the theorem.  \eproof

\mysection{Endomorphisms and their eigenvalues}{endomorphisms}

\myitpar{Representation of the polynomial algebra defined by an endomorphism.} Let $x$ be a variable, 
and let $\kk[x]$ be the $\kk$-algebra of polynomials in $x$ with coefficients in $\kk$\ffdot
Let $M$ be a $\kk$\dnsp-module\ffdot
The $\kk$\dnsp-endomorphisms $M\tto M$ form a $\kk$-algebra $\eend M$ with the composition as the multiplication.
For every $\kk$\dnsp-endomorphisms $E\colon M\tto M$ and every $a\elem\nn$ we will denote by $E^a$ 
the $a$\dnsp-fold composition $E\circ E\circ\ldots\circ E$\ffdot 
As usual, we interpret the $0$\dnsp-fold composition $E^0$ as the identity endomorphism $\id\elem\eend M$\ffdot
\newpage

For a $\kk$\dnsp-module endomorphism $E\colon M\tto M$ and a polynomial
\begin{equation*}
f(x)\off
=\off c_0 x^{n}  \qff +\qff c_1 x^{{n} {\minus}1} \qff +\qff c_2 x^{{n} {\minus}2} \qff +\qff \ldots \qff +\qff c_{n}
\dff\in\dff \kk[x]\endss
\end{equation*}
one can define an endomorphism $f(E)\colon M\tto M$ by the formula
\begin{equation*}
f(E)\qff 
=\qff c_0 E^{n}  \qff +\qff c_1 E^{{n} {\minus}1} \qff +\qff c_2 E^{{n} {\minus}2} \qff +\qff \ldots \qff +\qff c_{n}\endss.
\end{equation*}
The map $f(x)\mapsto f(E)$ is a homomorphism $\kk[x]\tto\eend M$ of $\kk$\dnsp-algebras.
This follows from the obvious identities $x^a\dff x^b\eeq x^{a{\plus}b}$ and $E^a\circ E^b\eeq E^{a{\plus}b}$\ffdot
This homomorphism defines a structure of $\kk[x]$\dnsp-module on $M$\ffdot 
Of course, this structure depends on $E$\ffdot

We will denote by $\kk[E]$ the image of the homomorphism $f(x)\mapsto f(E)$\dfdot
Since $\kk[x]$ is commutative, the image $\kk[E]$ is a commutative subalgebra of $\eend M$\dfdot

\myitpar{Eigenvalues.} Suppose that a $\kk$\dnsp-module endomorphism $E\colon M\tto M$ is fixed.

Let $\alpha\elem\kk$\ffdot
The kernel $\kker (E{\mminus}\alpha)$ is called the 
\emph{eigenmodule of\dss $E$ corresponding to} $\alpha$ and is denoted also by $E_{\alpha}$\ffdot
Clearly, $E_{\alpha}$ is a $\kk$\dnsp-submodule of $M$\ffdot
An element $\alpha\elem\kk$ is called an \emph{eigenvalue of\dss $E$} if the kernel $E_{\alpha}\eeq\kker (E{\mminus}\alpha)\nneq 0$\ffdot

The set of elements $v\elem M$ such than $(E{\mminus}\alpha)^i(v)\eeq 0$ for some $i\elem\nn$ is called the
\emph{extended eigenmodule of\dss $E$ corresponding to} $\alpha$ and is denoted by $\nnil(\alpha)$\ffdot
Clearly, $\nnil(\alpha)$ is a $\kk$\dnsp-submodule of $M$\ffdot

\mypar{Lemma.}{eigen} \emph{Let $\alpha\elem\kk$\ffdot 
Then the following statements hold.}
\vspace{-\bigskipamount}
\begin{itemize}
\item[\emph{\phantom{ii}(i)}] \hspace*{1.5em} \emph{The submodules $E_{\alpha}$ and $\nnil(\alpha)$ are $E$\dnsp-invariant.}
\vspace{-1ex} 
\item[\emph{\phantom{i}(ii)}] \hspace*{1.5em} \emph{$E_{\alpha}$ and $\nnil(\alpha)$ are $\kk[x]$\dnsp-submodules of $M$}\dfdot
\vspace{-1ex}
\item[\emph{(iii)}] \hspace*{1.5em} \emph{The submodule $\nnil(\alpha)$ is non-zero if and only if  $\alpha$ is an eigenvalue.}
\end{itemize}

\vspace{-\bigskipamount}
\proof {\fff}Let us prove (i)\hff, (ii)\sss first.
Note that $(E{\mminus}\alpha)^i\circ E =E\hff\circ (E{\mminus}\alpha)^i$ for every $i\elem\nn$\dnsp,
because $\kk[E]$ is a commutative subalgebra of $\eend M$\ffdot
Therefore, if $(E{\mminus}\alpha)^i(v)\eeq 0$\dfcom then
\[
(E{\mminus}\alpha)^i\left(E(v)\right)\off 
= \off\left((E{\mminus}\alpha)^i\circ E\right)\hff (v)\off 
=\off \left (E\hff\circ (E{\mminus}\alpha)^i\right) (v)\off 
=\off  E\hff \left((E{\mminus}\alpha)^i(v)\right)\off 
=\off 0\endss.
\]
In the case $i\eeq 1$ this implies that $E(E_{\alpha})\ssub E_{\alpha}$\ffdot
In general, this implies that 
\[
E(\kker (E{\mminus}\alpha)^i)\qff \ssub\qff \kker (E{\mminus}\alpha)^i\endss,
\]
and hence $E(\nnil (\alpha))\qff \ssub\qff \nnil(\alpha)$\ffdot
This proves (i)\hff, and (ii) immediately follows.\newpage

Finally,\qss let us prove (iii).\hspace*{0.4em} 
Suppose that $v\nneq 0$ and $(E{\mminus}\alpha)^i(v)\eeq 0$\dff\dfdot
Let $i$ be the smallest integer such that $(E{\mminus}\alpha)^i(v)\eeq 0$\dfdot
Note that $i\gres 0$ because $v\nneq 0$\ffdot
Let $w\eeq (E{\mminus}\alpha)^{i{\minus}1}(v)$\dnsp.\dss
Then $w\nneq 0$ and $(E{\mminus}\alpha)(w)\eeq 0$\dfdot
Therefore $E_{\alpha}\eeq\kker (E{\mminus}\alpha)\nneq 0$\dfdot
This proves (iii).  \eproof

\myitpar{Torsion free modules.} A $\kk$\dnsp-module $M$ is called \emph{torsion-free,}\sss
if $\alpha\dff m\eeq 0$ implies that either $\alpha\eeq 0$\dfcom or $m\eeq 0$\dfcom where $\alpha\elem\kk$ and $m\elem M$\dfdot
Since $\kk$ is assumed to be a ring without zero divisors, $\kk^n$ is a torsion free module for any non-zero $n\elem\nn$\dfdot

\vspace*{\bigskipamount}
\emph{\textbf{For the rest of this section we will assume that $M$ is a torsion-free module.}}

\mypar{Lemma.}{eigenmodules} \emph{Let\trs $\alpha_1\dff,\hspace{0.2em}\alpha_2\dff,\hspace{0.2em}\ldots\dff,\hspace{0.2em}\alpha_n$\dss
be distinct eigenvalues of an endomorphism\dss $E\colon M\tto M$\dnsp.\dss
Let\dss $Ker_1\dff,\hspace{0.2em}Ker_2\dff,\hspace{0.2em}\ldots\dff,\hspace{0.2em}Ker_n$\sss be the corresponding eigenmodules, i.e.
$Ker_i\qff =\qff \kker (E{\mminus}\alpha_i)$ 
for each $i\eeq 1\dff,\hspace{0.2em}2\dff,\hspace{0.2em}\ldots\dff,\hspace{0.2em}n$\nsp.\sss
Then the sum of these eigenmodules is a direct sum, i.e. an element
\[
v\dff\in\dff Ker_1\qff +\qff Ker_2\qff +\qff \ldots{\pplus}Ker_n 
\mbox{\ \ admits only one presentation\ \ }
v\off =\off v_1\qff +\qff v_2\qff +\qff \ldots\qff +\qff v_n
\]
with $v_i\dff\in\dff Ker_i$ for all $i\eeq 1\dff,\hspace{0.2em}2\dff,\hspace{0.2em}\ldots\dff,\hspace{0.2em}n$}\ffdot

\proof It is sufficient to prove that if $v_1{\pplus}v_2{\pplus}\ldots{\pplus}v_n\eeq 0$
and  $v_i\elem Ker_i$ for all $i$\ffcom
then $v_1\eeq v_2\eeq\ldots\eeq v_n\eeq 0$\ffdot
Suppose that $v_1{\pplus}v_2{\pplus}\ldots{\pplus}v_n\eeq 0$\ffcom $v_i\elem Ker_i$ for all $i$\ffcom and not all $v_i$ are equal to $0$\ffdot
Consider the maximal integer $m$ such that 
\begin{equation}
\label{min}
v_1\qff +\qff v_2\qff +\qff\ldots\qff +\qff v_m\eeq 0
\end{equation} 
for some elements $v_i\elem Ker_i$
such that $v_m\nneq 0$\ffdot
Note that in this case $v_i\nneq 0$ also for some $i\leqs m{\minus}1$\ffcom in view of (\ref{min})\fff.\sss
By applying $E{\mminus}\alpha_m$ to\sss (\ref{min})\fff,\dss we get
\[
(E{\mminus}\alpha_m)(v_1)\qff +\qff\ldots\qff +\qff
(E{\mminus}\alpha_m)(v_{m{\minus}1})\qff +\qff (E{\mminus}\alpha_m)(v_m)\off =\off 0\endss.
\]
Since $v_i\elem Ker_i\eeq\kker (E{\mminus}\alpha_i)$ and therefore $E(v_i)\eeq\alpha_i\fff v_i$ for all $i$\ffcom we see that
\begin{equation}
\label{assumed-min}
(\alpha_1{\mminus}\alpha_m)\fff v_1\qff +\qff\ldots\qff +\qff 
(\alpha_{m{\minus}1}{\mminus}\alpha_m)\fff v_{m{\minus}1}\qff +\qff
(\alpha_m{\mminus}\alpha_m)\fff v_m\off =\off 0\endss,
\end{equation}
\begin{equation}
\label{smaller}
(\alpha_1{\mminus}\alpha_m)\fff v_1\qff +\qff \ldots\qff +\qff
(\alpha_{m{\minus}1}{\mminus}\alpha_m)\fff v_{m{\minus}1}\off =\off 0\endss.
\end{equation}
Since the eigenvalues $\alpha_i$ are distinct, $\alpha_i{\mminus}\alpha_m\nneq 0$ for $i\leqs m{\minus}1$\ffdot
Since our module $M$ is assumed to be torsion-free, this implies that 
$(\alpha_i{\mminus}\alpha_m)\fff v_i\nneq 0$ if $i\leqs m{\minus}1$ and $v_i\nneq 0$\dnsp.\dss
As we noted above, $v_i\nneq 0$ for some $i\leqs m{\minus}1$\ffdot
Therefore, the equality (\ref{smaller}) contradicts to the choice of $m$\ffdot
This contradiction proves the lemma.  \eproof
\newpage

\mypar{Lemma.}{extended-eigenmodules} \emph{Let\dss $\alpha_1\dff,\hspace{0.2em}\alpha_2\dff,\hspace{0.2em}\ldots\dff,\hspace{0.2em}\alpha_n$\sss
be distinct eigenvalues of an endomorphism\dss $E\colon M\tto M$\dnsp.\dss
Let\dss $\nnnil_1\dff,\hspace{0.2em}\nnnil_2\dff,\hspace{0.2em}\ldots\dff,\hspace{0.2em}\nnnil_n$\sss be the corresponding extended eigenmodules, i.e.\sss
$\nnnil_i\eeq\nnil(\alpha_i)$\sss {\hspace*{0.2em}}for $i\eeq 1\dff,\hspace{0.2em}2\dff,\hspace{0.2em}\ldots\dff,\hspace{0.2em}n$\dnsp.\dss
Then the sum of these extended eigenmodules is a direct sum, i.e. an element
\[
v\dff\in\dff \nnnil_1\qff +\qff\nnnil_2\qff +\qff\ldots\qff +\qff\nnnil_n 
\mbox{\ \ admits only one presentation\ \ }
v\off =\off v_1\qff +\qff v_2\qff +\qff\ldots\qff +\qff v_n
\]
with $v_i\dff\in\dff \nnnil_i$ for all $i= 1\dff,\hspace{0.2em}2\dff,\hspace{0.2em}\ldots\dff,\hspace{0.2em}n$}\ffdot

\proof It is sufficient to prove that if $v_1{\pplus}v_2{\pplus}\ldots{\pplus}v_n\eeq 0$
and  $v_i\elem \nnnil_i$ for all $i$\ffcom
then $v_1\eeq v_2\eeq\ldots\eeq v_n\eeq 0$\ffdot
Suppose that $v_1{\pplus}v_2{\pplus}\ldots{\pplus}v_n\eeq 0$\ffcom $v_i\elem \nnnil_i$ for all $i$\ffcom and not all $v_i$ are equal to $0$\dfdot
The proof proceeds by replacing\fff,\dss in several steps (no more than $n$\dnsp)\hff,\sss the original elements $v_i$ by new ones in such a way 
that eventually not only $v_i\elem \nnnil_i$\ffcom  
but, moreover, $v_i\elem Ker_i\eeq\kker(E{\mminus}\alpha_i)$\dfcom
and still not all $v_i$ are equal to $0$\dfdot
Obviously, this will contradict to Lemma \ref{eigenmodules}.

Let $E_i\eeq E{\mminus}\alpha_i$ for all  $i\eeq 1\dff,\hspace{0.2em}2\dff,\hspace{0.2em}\ldots\dff,\hspace{0.2em}n$\dfdot
If $i\eeq 1\dff,\hspace{0.2em}2\dff,\hspace{0.2em}\ldots\dff,\hspace{0.2em}n{\minus}1$\dfcom or $n$\dfcom 
then $E_i^0(v_i)\nneq v_i$ and $E_i^{a}(v_i)\eeq 0$ for some integer $a\geqs 1$\dfdot 
If $v_i\nneq 0$\dfcom then we define $a_i$ as the largest integer $a\geqs 0$ such $E_i^{a}(v_i)\nneq 0$\dfdot  
Then $E_i^{a_i}(v_i)\nneq 0$ and  $E_i^{a_i{\plus}1}(v_i)\eeq 0$\dfdot
In particular, 
\[
(E{\mminus}\alpha_i)(E_i^{a_i}(v_i))\off
=\off E_i(E_i^{a_i}(v_i))\off
=\off E_i^{a_i{\plus}1}(v_i)\off
=\off 0\endss,
\]
and hence $E_i^{a_i}(v_i)\dff\in\dff Ker_i$\dfdot
If $v_i\eeq 0$\dfcom then we set $a_i\eeq 0$ and $E_i^{a_i}(v_i)\elem Ker_i$ is still true.

Let us fix an integer $k$ between $1$ and $n$\dfdot
Let $w_i\eeq E_k^{a_k}(v_i)$\dfcom where $i\eeq 1\dff,\hspace{0.2em}2\dff,\hspace{0.2em}\ldots\dff,\hspace{0.2em}n$\nsp.\sss
By applying $E_k^{a_k}$ to $v_1{\pplus}v_2{\pplus}\ldots{\pplus}v_n\eeq 0$\ffcom
we conclude $w_1{\pplus}w_2{\pplus}\ldots{\pplus}w_n = 0$\dfdot
Note that since the submodules $\nnnil_i$ are $E$\dnsp-invariant by Lemma \ref{eigen}, $w_i\elem\nnnil_i$ for every $i$\dfdot

{\sc Claim 1.}\dss \emph{If\dss $v_i\qff \neq\qff 0$\ffcom then $w_i\qff =\qff 0$\ffdot}

\emph{Proof of Claim 1.}\dss If $i\eeq k$ and $v_i\eeq v_k\nneq 0$\dfcom then $w_i\eeq w_k\eeq E_k^{a_k}(v_k)\nneq 0$  by the choice of $a_k$\dfdot
Suppose that $i\nneq k$ and $v_i\nneq 0$\dfdot 
Then 
\[
E_i^{a_i}(w_i)\off =\off E_i^{a_i}\left(E_k^{a_k}(v_i)\right)\off =\off E_k^{a_k}\left(E_i^{a_i}(v_i)\right)
\]
But $E_i^{a_i}(v_i)\elem Ker_i$ and  $E_i^{a_i}(v_i)\nneq 0$ by the choice of $a_i$\dfdot 
Since $E$ acts of $Ker_i$ as the multiplication by $\alpha_i$\dfcom we have
\[
E_k^{a_k}\left(E_i^{a_i}(v_i)\right)\off =\off (E{\mminus}\alpha_k)^{a_k}\left(E_i^{a_i}(v_i)\right)\off 
=\off (\alpha_i{\mminus}\alpha_k)^{a_k}\left(E_i^{a_i}(v_i)\right).
\]
Since $\alpha_i\nneq\alpha_k$ and $\kk$ is a ring without zero divisors, $(\alpha_i{\mminus}\alpha_k)^{a_k}\nneq 0$\dfdot
Since $M$ is a torsion free $k$\dnsp-module and $E_i^{a_i}(v_i)\nneq 0$\dfcom this implies that 
$(\alpha_i{\mminus}\alpha_k)^{a_k}\left(E_i^{a_i}(v_i)\right)\nneq 0$\dfdot
It follows that  
\[
E_i^{a_i}(w_i)\off =\off (\alpha_i{\mminus}\alpha_k)^{a_k}\left(E_i^{a_i}(v_i)\right)\off\neq\off 0\endss,
\]\newpage
and hence $w_i\nneq 0$\dfdot
This completes the proof of the claim.  \esubproof

{\sc Claim 2.}\dss \emph{If\dss $v_i\dff\in\dff Ker_i$\ffcom then  $w_i\dff\in\dff Ker_i$\ffdot}

\emph{Proof of Claim 2.}\dss Suppose that $v_i\elem Ker_i$\ffcom i.e. $E_i(v_i)\eeq 0$\dfdot
Since $E_i\eeq E{\mminus}\alpha_i$ and $E_k\eeq E{\mminus}\alpha_k$ obviously commute, 
it follows that
\[
E_i(w_i)\off =\off E_i(E_k^{a_k}(v_i))\off =\off  E_k^{a_k}(E_i(v_i))\off =\off 0\endss.
\]
This proves the claim.  \esubproof

To sum up, we see that by applying $E_k^{a_k}$ to the equality
$v_1{\pplus}v_2{\pplus}\ldots{\pplus}v_n\eeq 0$ with  $v_i\elem \nnnil_i$ for all $i$ 
we get another equality 
$w_1{\pplus}w_2{\pplus}\ldots{\pplus}w_n\eeq 0$ such that for all $i$\nsp:
\vspace{-2.4\bigskipamount}
\begin{itemize}
\item[{\phantom{ii}(i)}] \hspace*{1.5em} $w_i\elem \nnnil_i$\nsp; 
\vspace{-1ex}
\item[{\phantom{i}(ii)}] \hspace*{1.5em} \emph{if $v_i\nneq 0$\dfcom then $w_i\nneq 0$\nsp;}
\vspace{-1ex}
\item[{\phantom{}(iii)}] \hspace*{1.5em} \emph{if $v_i\elem Ker_i$\dfcom then $w_i\elem Ker_i$\nsp.} 
\end{itemize} 
\vspace{-\bigskipamount} 
In addition, $w_k\eeq E_k^{a_k}(v_k)\elem Ker_k$ even if $v_k$ did not belonged to the eigenmodule $Ker_k$\dnsp.\dss
Therefore, we can take $w_1\dff,\hspace{0.2em}w_2\dff,\hspace{0.2em}\ldots\dff,\hspace{0.2em}w_n$ as the new
elements $v_1\dff,\hspace{0.2em}v_2\dff,\hspace{0.2em}\ldots\dff,\hspace{0.2em}v_n$\nsp,\sss increasing the number of elements
belonging to the corresponding eigenmodules by an appropriate choice of $k$ (\fff if some $v_i$ did not belonged to eigenmodules yet\fff).

It follows that by starting with the equality $v_1{\pplus}v_2{\pplus}\ldots{\pplus}v_n\eeq 0$ and 
consecutively applying endomorphisms $E_k^{a_k}$ for $k\eeq 1\dff,\hspace{0.2em}2\dff,\hspace{0.2em}\ldots\dff,\hspace{0.2em}n$\ffcom
we will eventually prove the equality $v_1{\pplus}v_2{\pplus}\ldots{\pplus}v_n\eeq 0$ for some new vectors $v_i$
such that $v_i\elem Ker_i$ for all $i$\ffcom and still not all $v_i$ are equal to $0$\dfdot
The contradiction with Lemma \ref{eigenmodules} completes the proof.  \eproof

\mypar{Lemma.}{nil-basis} \emph{Let $E\colon M\tto M$ be an endomorphism of $M$ and let $\alpha$ be an eigenvalue of $E$\ffdot
Suppose that $v\elem\nnil(\alpha)$\dnsp.\sss 
Let $a\geqs 0$ be the largest integer such that $(E{\mminus}\alpha)^a(v)\nneq 0$\ffcom
and let $v_i\eeq (E{\mminus}\alpha)^i(v)$\dss for\dss $i\eeq 0\dff,\hspace{0.2em}1\dff,\hspace{0.2em}\ldots\dff,\hspace{0.2em}a$\ffdot
Then the homomorphism $\kk^{a{\plus}1}\tto M$ defined by
\[
(x_0\dff,\hspace{0.2em}x_1\dff,\hspace{0.2em}\ldots\dff,\hspace{0.2em}x_a)\off 
\longmapsto\off x_0 v_0 \qff +\qff x_1 v_1 \qff +\qff\ldots\qff +\qff x_a v_a
\]
is an isomorphism onto its image.
In particular, $v_0\dff,\hspace{0.2em}v_1\dff,\hspace{0.2em}\ldots\dff,\hspace{0.2em}v_a$ are free generators of a free $\kk$\dnsp-submodule of $M$\dfdot}

\proof It is sufficient to prove that our homomorphism is injective.
In other terms, it is sufficient to prove that if
\begin{equation}
\label{ldependence}
x_0 v_0 \qff +\qff x_1 v_1 \qff +\qff\ldots\qff +\qff x_a v_a\off =\off 0
\end{equation} 
for some $x_0\dff,\hspace{0.2em}x_1\dff,\hspace{0.2em}\ldots\dff,\hspace{0.2em}x_a\elem\kk$\ffcom
then $x_0\eeq x_1\eeq \ldots\eeq x_a\eeq 0$\ffdot
Suppose that (\ref{ldependence}) holds and $x_i\nneq 0$ for some $i$\ffdot
Let $b\elem\nn$ be the minimal integer with the property $x_b\nneq 0$\dfdot
Let us apply $(E{\mminus}\alpha)^{a{\minus}b}$ to (\ref{ldependence}).\dss
Note that if $i\gres b$\dfcom then 
\[
(E{\mminus}\alpha)^{a{\minus}b}(v_i)\off
=\off (E{\mminus}\alpha)^{a{\minus}b} \left( (E{\mminus}\alpha)^{i}(v) \right)\off
\neq\off 
(E{\mminus}\alpha)^{a{\minus}b{\plus}i}(v)\off
=\off 0
\]
because ${a{\minus}b{\plus}i}\gres a$ and $(E{\mminus}\alpha)^n(v)\eeq 0$ for $n\gres a$ by the choice of $a$\dfdot
Therefore, the operator $(E{\mminus}\alpha)^{a{\minus}b}$ takes the left hand side of (\ref{ldependence}) to
\[
x_b (E{\mminus}\alpha)^{a{\minus}b}(v_b)\off 
=\off x_b (E{\mminus}\alpha)^{a{\minus}b}\left( E{\mminus}\alpha)^{b}(v)\right)\off 
=\off x_b (E{\mminus}\alpha)^{a{\minus}b{\plus}b}(v)\off 
=\off x_b (E{\mminus}\alpha)^{a}(v)\endss,
\]
and hence the result of application of $(E{\mminus}\alpha)^{a{\minus}b}$ to (\ref{ldependence}) is
\begin{equation}
\label{one-term}
x_b (E{\mminus}\alpha)^{a}(v)\off =\off 0\endss.
\end{equation}
But $(E{\mminus}\alpha)^{a}(v)\nneq 0$ by the choice of $a$\dfcom
and $x_b\nneq 0$ by the choice of $b$\dfdot
Since the module $M$ is assumed to be torsion free, these facts together with (\ref{one-term}) lead to a contradiction.
This contradictions shows that (\ref{ldependence}) may be true only if $x_i\eeq 0$ for all $i$\ffdot  \eproof

\mysection{Torsion modules and a property of free modules}{free-modules}

\myitpar{Torsion modules.} An element $m\elem M$ of a $\kk$\dnsp-module $M$ is called 
a \emph{torsion element}\dss if $x\dff m\eeq 0$ for some non-zero $x\elem\kk$\dfdot
A $\kk$\dnsp-module $M$ is called a \emph{torsion module} if every element of $M$ is a torsion element.

\mypar{Lemma.}{entire-rings} \emph{Let $n\elem\nn$\dfdot
If\dss $M$ is a $\kk$\dnsp-submodule of a $\kk$\dnsp-module $N$ and both $M$ and $N$ are isomorphic to $\kk^n$\dfcom 
then the quotient $N/M$ is a torsion module.}

\proof Suppose that $N/M$ is not a torsion module. 
Then there is an element $v\elem N/M$ such that $\alpha\dff v\nneq 0$ if $\alpha\nneq 0$\dfdot
For such a $v$ the map $\alpha\mapsto \alpha\dff v$ is an injective homomorphism of $\kk$\dnsp-modules $\kk\tto N/M$\dfdot
Let us lift $v\elem N/M$ to an element $v_0\elem N$\dfcom
so $v$ is the image of $v_0$ under the canonical surjection $N\tto N/M$\dfdot
Then the map $\alpha\mapsto \alpha\dff v_0$ is an injective homomorphism of $\kk$\dnsp-modules $\kk\tto N$\dfdot

Clearly, if $v_1\dff,\dff v_2\dff,\dff\ldots\dff,\dff v_n$ is a basis of $M$ (which exists because $M$ is isomorphic to $\kk^n$\nsp),\dss
then $v_0\dff,\dff v_1\dff,\dff\ldots\dff,\dff v_n$ is a basis of $\kk v_0 {\pplus}M$\dfdot
Therefore, $\kk v_0 {\pplus}M$ is a submodule isomorphic to $\kk^{n{\plus}1}$ of the module $N$ isomorphic to $\kk^n$\dfdot
In particular, there exist an injective $\kk$\dnsp-homomorphism $J\colon\kk^{n{\plus}1}\tto\kk^n$\dfdot

Since $\kk$ has no zero divisors, it can be embedded into its field of fractions, which we will denote by $\field$\dfdot
Moreover, the $\kk$\dnsp-homomorphism $\kk^{n{\plus}1}\tto\kk^n$ extends to an $\field$\dnsp-linear map
$\field^{n{\plus}1}\tto\field^n$\dfcom which we will denote by $J_F$\dfdot

{\sc Claim.} $J_F$ \emph{is injective.}

\emph{Proof of the claim.} Suppose $(y_0\dff,\dff y_1\dff,\dff\ldots\dff,\dff y_n)\elem\field^{n{\plus}1}$ 
is non-zero and belongs to the kernel of $J_F$\dfdot
Since $\field$ is the field of fractions of $\kk$\dfcom
there is an element $z\elem\kk$ such that $z\dff y_0\dff,\dff z\dff y_1\dff,\dff\ldots\dff,\dff z\dff y_n\elem\kk$\dfdot
For such an element $z\elem\kk$ the $(n{\plus}1)$\dnsp-tuple $(z\dff y_0\dff,\dff z\dff y_1\dff,\dff\ldots\dff,\dff z\dff y_n)$ 
belongs to $\kk^{n{\plus}1}$\dfcom and 
\begin{align*}
J(z\dff y_0\dff,\dff z\dff y_1\dff,\dff\ldots\dff,\dff z\dff y_n)\off 
&=\off J_F (z\dff y_0\dff,\dff z\dff y_1\dff,\dff\ldots\dff,\dff z\dff y_n)\\
&=\off z\trf J_F (y_0\dff,\dff y_1\dff,\dff\ldots\dff,\dff y_n)\off =\off z\dff 0\off =\off 0\endss.
\end{align*}
Since $\field$ is the field of fractions of $\kk$\dfcom
$(y_0\dff,\dff y_1\dff,\dff\ldots\dff,\dff y_n)\nneq 0$  implies that the $(n{\plus}1)$\dnsp-tuple 
$(z\dff y_0\dff,\dff z\dff y_1\dff,\dff\ldots\dff,\dff z\dff y_n)\nneq 0$\dfdot
At the same time this $(n{\plus}1)$\dnsp-tuple belongs to the kernel of $J$\dfcom
in contradiction with the injectivity of $J$\dfdot
The claim follows.  \esubproof

As is well known, for a field $\field$ there are no injective $\field$\dnsp-linear maps $\field^{n{\plus}1}\tto\field^n$\dfdot 
The contradiction with the above claim proves that $N/M$ is indeed a torsion module.  \eproof

\mysection{Polynomials and their roots}{polynomials}

\mypar{Lemma.}{roots} \emph{Let $p(x)\elem\kk[x]$ be a polynomial with leading coefficient $1$\dfcom and let $\alpha\elem\kk$\dfdot
Then $\alpha\elem\kk$ is a root of $p(x)$ if and only if} 
\begin{equation}
\label{linear-factor}
p(x)\off =\off (x-\alpha)q(x)
\end{equation}
\emph{for some polynomial $q(x)\elem\kk[x]$ with leading coefficient $1$\dfdot
If $\alpha$ is a root, then $q(x)$ is uniquely determined by\dss \textup{(\ref{linear-factor})}.}

\proof Suppose that $p(x)\eeq (x-\alpha)q(x)$ and both $p(x)\dff,\dff q(x)$ have the leading coefficients $1$\dfdot
Then $\deg q(x)\eeq \deg p(x){\minus}1$ and the polynomials $p(x)\dff,\dff q(x)$ have the form
\begin{align*}
p(x)\off 
&=\off x^n\qff +\qff c_1\fff x^{n{\minus}1\fff}\qff +\qff \ldots\qff +\qff a_{n{\minus}1}\fff x\qff +\qff c_{n}\endss,\\
q(x)\off 
&=\off x^{n{\minus}1}\qff +\qff d_1\fff x^{n{\minus}2}\qff +\qff \ldots\qff +\qff d_{n{\minus}2}\fff x\qff +\qff d_{n{\minus}1}\endss,
\end{align*}
where $c_1\dff,\dff\ldots\dff,\dff c_n\dff,\dff d_1\dff,\dff\ldots\dff,\dff d_{n{\minus}1}\elem\kk$\dfdot
Let us compute the product 
\[
(x-\alpha)q(x)\off 
=\off (x-\alpha) \bigl(x^{n{\minus}1}\qff +\qff \ldots\qff +\qff \phantom{\alpha\fff}d_{n{\minus}2}\fff x\phantom{{\fff\hff}^2}\qff +\qff \phantom{\alpha\fff}d_{n{\minus}1}\fff\bigr).
\]

Obviously\nsp,\qss
\begin{align*}
(x-\alpha)q(x)\off 
               &=\off  \phantom{(}x^{n}\qff +\qff d_1\fff x^{n{\minus}1}\qff +\qff \ldots\qff +\qff \phantom{\alpha\fff}d_{n{\minus}2}\fff x^2\qff +\qff \phantom{\alpha\fff}d_{n{\minus}1}\fff x \\
               &\off\phantom{=  (x^{n}\qff +}\nsp-\qff \alpha\fff x^{n{\minus}1}\qff  -\qff \ldots\qff -\qff \alpha\fff d_{n{\minus}3}\fff x^2\qff -\qff \alpha\fff d_{n{\minus}2}\fff x\qff -\qff  \alpha\fff d_{n{\minus}1}\endss. 
\end{align*} 
It follows that $p(x)\qff =\qff (x-\alpha)q(x)$ if and only if
\begin{align*}
c_1\off            &=\off d_1\qff -\qff \alpha \\
c_2\off            &=\off d_2\qff -\qff \alpha\fff d_1 \\
\ldots\off         &=\off \ldots\qff\ldots \\
c_{n{\minus}1}\off &=\off d_{n{\minus}1}\qff -\qff \alpha\fff d_{n{\minus}2} \\
c_n\off            &=\off -\qff\alpha\fff d_{n{\minus}1} \endss,
\end{align*}
or, equivalently,
\begin{align*}
d_1\off              &=\off \alpha\qff +\qff c_1 \\
d_2\off              &=\off \alpha\fff d_1\qff +\qff c_2 \\
\ldots\off           &=\off \ldots\qff\ldots \\
d_{n{\minus}1}\off   &=\off \alpha\fff d_{n{\minus}2}\qff +\qff c_{n{\minus}1}  \\
\alpha\fff d_{n{\minus}1}\qff +\qff c_n\off &=\off 0 \endss.
\end{align*}
These equalities allow to compute the coefficients $d_1\dff,\dff d_2\dff,\dff\ldots\dff,\dff d_{n{\minus}1}$
in terms of the coefficients $c_1\dff,\dff c_2\dff,\dff\ldots\dff,\dff c_{n{\minus}1}$\dnsp.\dss
Namely, $d_1\qff =\qff \alpha + c_1 $ and
\begin{equation*}
d_{\fff i}\off 
=\off \alpha^i\qff +\qff c_{1}\fff\alpha^{i{\minus}1}\qff  
+\qff c_{2}\fff\alpha^{i{\minus}2}\qff +\qff \ldots\qff +\qff c_{i{\minus}1}\fff\alpha\qff +\qff c_i
\end{equation*}
for $2\leqs i\leqs n{\minus}1$\dfdot
Therefore, the last equality $\alpha\fff d_{n{\minus}1}{\pplus}c_n\eeq 0$ holds if and only if
\[
\alpha\bigl(\alpha^{n{\minus}1}\qff +\qff c_1\fff\alpha^{n{\minus}2}\qff 
+\qff \ldots\qff c_{n{\minus}1}\fff\alpha\bigr)\qff + \qff c_n\off =\off 0\endss,
\]
i.e. if and only if $p(\alpha)\qff =\qff 0$\dfdot
The lemma follows.  \eproof

\mypar{Corollary.}{mult-roots} \emph{Let\dss $p(x)\elem\kk[x]$\dss be a polynomial with leading coefficient\dss $1$\dfcom 
and let\dss $\alpha\elem\kk$\dnsp.\dss
Then there is a number\dss $m\elem\nn$\dss and a polynomial\trs $r(x)\elem\kk[x]$\dss such that\dss 
$p(x)\eeq (x{\mminus}\alpha)^m r(x)$\dss and\dss $\alpha$\dss is not a root of\dss $r(x)$\dfdot
The number\dss $m$\dss and the polynomial\dss $r(x)$\dss are uniquely determined by\dss $p(x)$\dss and\dss $\alpha$\dfdot}

\proof If $p(\alpha)\nneq 0$\dfcom then, obviously, $m\eeq 0$ and $r(x)\eeq p(x)$\dfdot
If $p(\alpha)\eeq 0$\dfcom we can apply Lemma \ref{roots}.
If $q(\alpha)\nneq 0$\dfcom then $m\eeq 1\dff,\dff r(x)\eeq q(x)$ and we are done.
If $q(\alpha)\eeq 0$\dfcom then we can apply Lemma \ref{roots} again.
Eventually we will get a presentation $p(x)\eeq (x{\mminus}\alpha)^m r(x)$ such that $r(\alpha)\nneq 0$\dfdot
By consecutively applying the uniqueness part of Lemma \ref{roots}, 
we see that $(x{\mminus}\alpha)^{m{\minus}1} r(x)$\dfcom
$(x{\mminus}\alpha)^{m{\minus}2} r(x)\dff,$ \ldots\dff,\dss
and, eventually, $m$ and $r(x)$ are uniquely determined by $p(x)$ and $\alpha$\dfdot  \eproof

\myitpar{The multiplicity of a root.} We will denote by $\deg\fff p(x)$ the degree of the polynomial $p(x)$\dfdot
If $\alpha$ is a root of $p(x)$\dfcom 
then the number $m$ from Corollary \ref{mult-roots} is called the\dss \emph{multiplicity}\dss of the root $\alpha$\dfdot

\mypar{Corollary.}{root-factors} \emph{Let\dss $p(x)\elem\kk[x]$\dss be a polynomial with leading coefficient\dss $1$\dfdot
The number\dss $k$\dss of distinct roots of\dss $p(x)$\dss is finite and\dss $k\leqs \deg p(x)$\dfdot
If\dss $\alpha_1\dff,\dff\alpha_2\dff,\dff\ldots\dff,\dff\alpha_k$\dss is the list of all distinct roots of\dss $p(x)$\dfcom
and if\dss $\mu_1\dff,\dff\mu_2\dff,\dff\ldots\dff,\dff\mu_k$\dss are the respective multiplicities of these roots,
then}
\begin{equation}
\label{full-factorization}
p(x)\off 
=\off \bigl(x-\alpha_1\bigr)^{\mu_1} \bigl(x-\alpha_2\bigr)^{\mu_2} \cdots \bigl(x-\alpha_k\bigr)^{\mu_k}\dff r(x)\endss,
\end{equation}
\emph{where\dss $r(x)\elem\kk[x]$\dss has no roots in\dss $\kk$\dfdot
The polynomial\dss $r(x)$\dss is uniquely determined by $p(x)$\dfdot}

\proof By consecutively applying the existence part of the Corollary \ref{mult-roots},
we see that a factorization of the form (\ref{full-factorization}) exists. 
Similarly, the uniqueness of $r(x)$ follows from the uniqueness part of Corollary \ref{mult-roots}.\eproof

\myitpar{Polynomials with all roots in\dss $\kk$\dfdot} Again, let $p(x)\elem\kk[x]$\dss 
be a polynomial with leading coefficient\dss $1$\dfdot
We say that $p(x)$\sss \emph{has all roots in}\dss $\kk$\dfcom
if in the factorization (\ref{full-factorization}) the polynomial $r(x)\eeq 1$\dfdot
In other words, $p(x)$\sss \emph{has all roots}\dss in $\kk$\dfcom if $p(x)$ has the form
\[
p(x)\off 
=\off \bigl(x-\alpha_1\bigr)^{\mu_1} \bigl(x-\alpha_2\bigr)^{\mu_2} \cdots \bigl(x-\alpha_k\bigr)^{\mu_k}\endss,
\]
for some $\alpha_1\dff,\dff\alpha_2\dff,\dff\ldots\dff,\dff\alpha_k\dff\elem\kk$ and some non-zero 
$\mu_1\dff,\dff\mu_2\dff,\dff\ldots\dff,\dff\mu_k\dff\elem\nn$\dfdot
Obviously, then
\[
\deg\fff p(x)\off =\off \mu_1{\pplus}\mu_2{\pplus}\ldots{\pplus}\mu_k\endss.
\]

\mysection{The main theorems}{relations}

Let us fix for the rest of this section a polynomial 
\[
p(x)\off =\off x^n\qff +\qff c_1\fff x^{n{\minus}1}\qff +\qff \ldots\qff +\qff 
c_{n{\minus}1}\fff x\qff +\qff c_{n}\dff\in\dff\kk[x]
\]
with the leading coefficient\sss $1$\dfdot

Consider the left shift operator $L\colon\seqq_\kk\tto\seqq_\kk$ from Section \ref{seq}.
As it was explained in Section \ref{endomorphisms}, the operator $L$ defines a homomorphism of $\kk$-algebras $\kk[x]\tto\eend\seqq_\kk$
by the rule $f(x)\mapsto f(L)$\dfdot
We are interested in the kernel $\kker p(L)$\dfdot
\newpage

\mypar{Lemma.}{kernel-recurr} \emph{A sequence\sss $s\elem\seqs$\sss belongs to\dss $\kker p(L)$\dss if and only if
\begin{equation}
\label{classical-recurr}
s_{\fff i}\qff +\qff c_1\fff s_{\fff i{\minus}1\fff}\qff 
+\qff \ldots\qff +\qff c_{n{\minus}1}\fff s_{\fff i{\minus}n{\plus}1}\qff +\qff c_{n}\fff  s_{\fff i{\minus}n}\off 
=\off 0
\end{equation}
for all\dss $i\elem\nn$\dfcom {\dff}$i\geqs n$\dfdot}

\proof Let us compute the terms of $p(L)(s)$\dfcom using Corollary \ref{iterated-shift} at the last step:
\begin{align*}
\bigl(\fff p(L) (s)\bigr)_i\off &=\off \dff 
\bigl(\dff L^n\qff +\qff c_1\hff L^{ n {\minus} 1}\qff +\qff \ldots\qff 
+\qff c_{n{\minus}1}\hff L\qff +\qff c_n )\hff (s)\dff\bigr)_{i}\vspace*{\medskipamount} \\
                            &=\off \dff 
                            \bigl(L^n (s) \bigr)_{i}\qff +\qff c_1\fff \bigl(L^{ n {\minus} 1} (s) \bigr)_{i}\qff 
                            +\qff  \ldots\qff +\qff
                               c_{n{\minus}1}\fff\bigl(L (s) \bigr)_{i}\qff +\qff  c_{n}\fff  \bigl( s \bigr)_{i}\vspace*{\medskipamount} \\
                            &=\off \hspace*{0.3em}s_{\fff i{\plus}{n}}\qff +\qff 
                            c_1\fff s_{\fff i{\plus}{n}{\minus}1}\qff +\qff \ldots\qff +\qff 
                               c_{n{\minus}1}\fff s_{\fff i{\plus}1}\qff +\qff c_{n}\fff  s_{\fff i}\endss.\vspace*{\medskipamount} 
\end{align*}
This calculation shows that $s\elem\kker f(L)$ if and only if 
\begin{equation}
\label{shift-classical-recurr}
s_{\fff i{\plus}{n}}\qff +\qff c_1\fff s_{\fff i{\plus}{n}{\minus}1}\qff +\qff \ldots\qff +\qff 
c_{n{\minus}1}\fff s_{\fff i{\plus}1}\qff +\qff c_{n}\fff  s_{\fff i}\off =\off 0
\end{equation}
for all integers $i\geqs 0$\dfdot
Clearly,\dss (\ref{shift-classical-recurr})\dss holds for all integers $i\geqs 0$ if and only if\dss 
(\ref{classical-recurr})\dss holds for all integers $i\geqs n$\dnsp.\dss
The lemma follows.  \eproof

\myitpar{Remark.} Classically, a sequence $s\elem\seqs$ is called\dss \emph{recurrent}\trs
if its terms satisfy (\ref{classical-recurr}) for all $i\elem\nn$\dfcom $i\geqs n$\dfcom
and the equation\sss (\ref{classical-recurr})\sss is called a\dss \emph{linear recurrence relation}.\dss
This explains the title of the paper.

\mypar{Lemma.}{freeness} \emph{The map\dss $F\colon\kker p(L)\tto\kk^n$\dss defined by}
\[
F\colon s\qff\longmapsto\qff (s_0\dff,\dff s_1\dff,\dff \ldots\dff,\dff s_{n{\minus}1})\dff\in\dff k^n
\] 
\emph{is an isomorphism. 
In particular,\dss $\kker p(L)$\dss is a free\dss $\kk$\dnsp-module of rank\dss $n$\dfdot}

\proof By Lemma \ref{kernel-recurr} the kernel $\kker p(L)$ is equal to the $\kk$\dnsp-submodule of $\seqs$ consisting of sequences $s$
satisfying the relation (\ref{classical-recurr}) for all\dss $i\elem\nn$\dfcom {\dff}$i\geqs n$\dfdot
Clearly, (\ref{classical-recurr}) allows to compute each term $s_i\dff,\dff i\geqs n$ of $s$ as the linear combination
of $n$ immediately preceding terms $s_{i{\minus}1\fff}\dff,\dff s_{i{\minus}2\fff}\dff,\dff \ldots\dff,\dff s_{i{\minus}n\fff}$ of $s$
with coefficients $c_1\dff,\dff c_2\dff,\dff \ldots\dff,\dff c_n$ independent of $i$\dfdot
Therefore, such a sequence $s$ is determined by its first $n$ terms $s_0\dff,\dff s_1\dff,\dff \ldots\dff,\dff s_{n{\minus}1}$\dfdot
Moreover, these $n$ terms can be prescribed arbitrarily.
The lemma follows.  \eproof

\mypar{Theorem.}{powers-of-roots} \emph{Suppose that\dss $\alpha\elem\kk$\dss 
is a root of $p(x)$ of multiplicity\dss $\mu$\dnsp.\dss
Then\dss $p(L)(s(\alpha\dff,\dff a))\eeq 0$\dss for each $a\elem\nn$\dfcom $0\leqs a\leqs\mu{\minus}1$\dfcom
where\dss $s(\alpha\dff,\dff a)$\dss are the sequences defined in
Section \ref{operators}, the paragraph immediately preceding Lemma \ref{eigen-powers}.}

\proof By Corollary \ref{mult-roots}\hff, $p(x)$ has the form
$p(x) = (x-\alpha)^\mu\dff q(x)$\dfdot
Therefore
\begin{equation}
\label{p0}
p(L)\off =\off \bigl(L\qff -\qff \alpha\bigr)^\mu\dff q(L)\off =\off q(L)\dff \bigl(L\qff -\qff\alpha\bigr)^\mu\endss.
\end{equation}
If $a\leqs\mu{\minus}1$\dfcom then $\mu{\minus}1{\minus}a\geqs 0$\dfcom 
and we can present $\bigl(L-\alpha\bigr)^\mu$ as the following product:
\begin{equation}
\label{p1}
\bigl(L\qff -\qff \alpha\bigr)^\mu\off 
= \off\bigl(L-\alpha\bigr)^{\mu{\minus}1{\minus}a} \bigl(L\qff -\qff \alpha\bigr) 
\bigl(L\qff -\qff \alpha\bigr)^a\endss.
\end{equation}
By Lemma \ref{left-der}
\begin{equation}
\label{p2}
\bigl(L\qff -\qff\alpha\bigr)^a\dff \bigl(s(\alpha\dff,\dff a)\bigr)\off 
= \off s(\alpha\dff,\dff a{\mminus}a)\off 
=\off s(\alpha\dff,\dff 0)\off 
=\off s(\alpha)\endss,
\end{equation}
and by Lemma \ref{eigen-powers}, 
\begin{equation}
\label{p3}
\bigl(L\qff -\qff\alpha\bigr)\bigl(s(\alpha)\bigr)\off 
=\off 0\endss.
\end{equation}
By combining\dss (\ref{p1}),\dss (\ref{p2}),\dss and\dss (\ref{p3}),\dss we see that
\begin{equation*}
\bigl(L\qff -\qff\alpha\bigr)^\mu\dff \bigl(s(\alpha\dff,\dff a)\bigr)\off 
=\off 0\endss.
\end{equation*}
By combining the last equality with\dss (\ref{p0}),\dss we get
\begin{equation*}
p(L) \bigl(s(\alpha\dff,\dff a)\bigr)\off 
=\off q(L) \bigl(L\qff -\qff \alpha\bigr)^\mu\dff \bigl(s(\alpha\dff,\dff a)\bigr)\off 
=\off q(L) (0)\off 
=\off 0\endss.
\end{equation*}
The theorem follows.  \eproof

\mypar{Theorem.}{all-roots-solutions} \emph{Suppose that\dss $p(x)$\dss has all roots in $\kk$\dfdot 
Let\dss $k$\dss be the number of distinct roots of\dss $p(x)$\dfcom
let\dss $\alpha_1\dff,\dff\alpha_2\dff,\dff\ldots\dff,\dff\alpha_k$\dss be these roots, 
and let \dss $\mu_1\dff,\dff\mu_2\dff,\dff\ldots\dff,\dff\mu_k$\dss be, respectively\fff,\sss the multiplicities of these roots.
Then sequences\dss $s(\alpha_u, a)$\dfcom where\dss $1\leqs u\leqs k$\dss and\dss $0\leq a\leq \mu_u{\minus}1$\dfcom
are free generators of a free\sss $\kk$\dnsp-submodule of\dss $\kker f(L)$\dnsp.}

\proof By Theorem \ref{powers-of-roots}, all these sequences belong to $\kker p(L)$\dnsp.\dss
In particular, they are generators of a $\kk$\dnsp-submodule of $\kker p(L)$\dfdot
Let us prove that they are free generators.

Let $1\leqs u\leqs k$\dfdot
By Lemma \ref{left-der}, 
\begin{equation}
\label{down}
\bigl( L{\mminus}\alpha_u \bigr)^{\mu_u{\minus}1}\fff\bigl(  s(\alpha_u\dff,\dff \mu_u{\minus}1)\bigr)\off 
=\off s(\alpha_u\dff,\dff (\mu_u{\minus}1){\minus}(\mu_u{\minus}1))\off
=\off s(\alpha_u\dff,\dff 0) = s(\alpha_u)\endss.
\end{equation}
By Lemma \ref{eigen-powers}, $\bigl(L{\mminus}\alpha_u\bigr)\bigl(s(\alpha_u)\bigr)\eeq 0$ and hence
\begin{equation*}
(L{\mminus}\alpha_u)^{\mu_u}\dff s(\alpha_u\dff,\dff \mu_u{\minus}1)\off 
=\off (L{\mminus}\alpha_u)\dff(s(\alpha))\off 
=\off 0\endss.
\end{equation*}
In particular,\dss $s(\alpha_u\dff,\dff \mu_u{\minus}1)$\sss belongs 
to the extended eigenmodule\dss $\nnil(\alpha_u)$\dfdot  
\newpage

In addition, (\ref{down}) together with the fact that $s(\alpha_u)\nneq 0$ 
implies that $\mu_u{\minus}1$ is the largest integer $a$ such that 
\[
\bigl(L{\mminus}\alpha_u\bigr)^a\fff\bigl( s(\alpha_u\dff,\dff \mu_u{\minus}1)\bigr)\off
\neq\off 0\endss.
\]
Lemma \ref{nil-basis} implies that the sequences $(L{\mminus}\alpha_u)^a\fff\bigl( s(\alpha\dff,\dff \mu_u{\minus}1)\bigr)$ for 
$a\eeq 0\dff,\hspace{0.2em}1\dff,\hspace{0.2em}\ldots\hspace{0.2em},\dff\mu_u{\minus}1$
form a basis of a free submodule of $\nnil(\alpha_u)\ssub\kker p(L)$\dfdot
Since
\[
\bigl(L{\mminus}\alpha_u\bigr)^a\bigl(s(\alpha_u\dff,\dff \mu_u)\bigr)\off
=\off s(\alpha_u\dff,\dff \mu_u{\minus}a)
\]
by Lemma \ref{left-der}, this implies that the sequences\sss $s(\alpha_u\dff,\dff a)$\sss for\sss 
$a\eeq 0\dff,\hspace{0.2em}1\dff,\hspace{0.2em}\ldots\hspace{0.2em},\dff\mu_u{\minus}1$\sss
form a basis of a free submodule of\dss $\nnil(\alpha_u)\ssub\kker p(L)$\dfdot
By combining this result with Lemma \ref{extended-eigenmodules}, we see that
the sequences\sss $s(\alpha_u\dff,\dff a)$\sss from the theorem form a basis of a free submodule of $\kker p(L)$\dfdot
This completes the proof.  \eproof

\mypar{Theorem.}{torsion-quot-delta} \emph{Let\dss $S\ssub\nsp\kker p(L)$\dss be the free $\kk$\dnsp-module generated by 
the sequences\sss $s(\alpha_u\dff,\dff a)$\sss from Theorem \ref{all-roots-solutions}.
Then the quotient\sss $\kk$\dnsp-module\dss $(\kker p(L))/S$\sss is a torsion module.}

\proof Let $n\eeq\deg\fff p(x)$\dfdot
By Lemma \ref{freeness}, $\kker p(L)$ is a free module of rank $n$\dfcom
i.e. is isomorphic to $\kk^n$\dfdot
Since $n\eeq\mu_1{\pplus}\ldots{\pplus}\mu_k$\dfcom
we have exactly $n$ sequences\sss $s(\alpha_u\dff,\dff a)$\dnsp.\dss
By Theorem \ref{all-roots-solutions}, they are free generators of $S$\dfdot
In particular, $S$ is also isomorphic to $\kk^n$\dfdot
It remains to apply Lemma \ref{entire-rings}.  \eproof

\mypar{Corollary.}{vector-spaces} \emph{If\dss $\kk$\sss is a field, then $S\qff =\qff \nsp\kker p(L)$\dfdot}

\proof A torsion module over a field is equal to $0$\dfdot  \eproof

\mynonumsection{Note bibliographique}

\vspace*{\bigskipamount}
This note is concerned only with the works which influenced the present paper.\qss
The author did not attempted to write even an incomplete account of
the history of the theory of linear recurrence relations.

The theory of recurrent sequences is an obligatory topic for any introduction to combinatorics.\qss
But all too often the proofs are presented only for Fibonacci numbers,\qss 
even if the general case is discussed.\qss
The author stumbled upon this tradition in P. Cameron's textbook\qss \cite{c}.\qss
The discussion of the general case in\qss \cite{c} is limited by the following. 
\vspace*{-\bigskipamount}

\begin{quote}
In this case,\qss suppose that $\alpha$ is a root of the characteristic equation with multiplicity $d$\dfdot
Then it can be verified that the the $d$ functions 
$\alpha^n$\dnsp,\qss $n\dff \alpha^n$\dnsp,\qss $\ldots$,\qss $n^{d{\minus}1}\dff\alpha^n$
are all solutions of of the recurrence relation.\qss
Doing this for every root,\qss
we again find enough independent solutions that $k$ initial values can be fitted.\qss

The justification of this is the fact that the solutions claimed 
can be substituted in the recurrence relation and its truth verified.
\end{quote}

\vspace*{-\bigskipamount}
This discussion ignores at least two significant issues.\qss
First,\qss the claim that the sequences $n^{\hff i}\dff\alpha^n$ provide solutions
cannot be justified by the substitution of them in the recurrence relations and a routine verification
simply because they are solutions only if $i\leq d{\minus}1$\dfdot 
Second,\qss one needs to prove that these solutions are linearly independent.\qss

The standard approach to the general case is based on the theory of generating functions
and the partial fractions expansion of rational functions.\qss
This method is elegantly presented in Chapter 3 of M. Hall's classical book \cite{hall}.\qss
A recent presentation of this method can be found in Section 3.1 of M. Aigner's book \cite{a}.\qss
When this approach is chosen,\qss the existence of partial fraction expansions is simply quoted
as a tool external to the theory of linear recurrence relations.\qss

The theory of partial fractions is an application of the
theory of modules over principal entire rings.\qss
See N. Bourbaki\qss \cite{b}\fff,\qss Section 2.{\fff}3,\qss or\qss 
S. Lang\qss \cite{l}\fff,\qss Section IV.{\fff}5{\fff}.\qss
Hence the theory of linear recurrence relations is also an application of the theory of modules.\qss
Once this is realized,\qss it is only natural to use the theory of modules directly,\qss
without using the generating functions and partial fractions as an intermediary.\qss

Another application of the same part of the theory of modules is the theory of the Jordan normal form.\qss
See N. Bourbaki\qss \cite{b}\fff,\qss Section 5,\qss or\qss S. Lang,\qss Chapter XIV.\qss
Of course,\qss the theory of the Jordan normal form precedes the theory of modules
and is usually presented without any references to the latter.\qss

It is only natural to adapt directly the standard arguments 
from the theory of the Jordan normal form
to the theory of linear recurrence relation.\qss
This is done in Section\qss \ref{endomorphisms}\qss of the present paper{},\qss
which was heavily influenced by I.M. Gelfand's classics \cite{gelfand}.\qss

The material of Sections\dss \ref{seq}\dss and\dss \ref{operators}\dss took its 
present form under the influence of G.-C. Rota's ideas about the umbral calculus.\qss
See,\qss for example,\qss \cite{rr}\fff.

The definition of divided derivatives was motivated by the desire to
prove the main results without any restrictions on the characteristic of the base ring\qss $\kk$\qss and
at the same time to avoide brute force calculations with binomial coefficients.\footnote{For the purposes
of this paper{},\qss any calculation with binomial coefficients using the well known expression
in terms of factorials is a brute force calculation.}
After this work was completed,\qss the author came across 
the paper\qss \cite{d}\qss by J. Dieudonn\'{e},\qss
from which he learned that this notion was first introduced by 
H. Hasse,\qss F.K. Schmidt,\qss and O. Teichm\"{u}ller in 1936,\qss 
with applications to various questions of algebra in mind.\qss
See the references in\qss \cite{d}.\qss
This was a pleasant surprise,\qss especially because a major part of the author's work
is devoted to Teichm\"{u}ller modular groups and Teichm\"{u}ller spaces.

\vspace*{1ex}

\begin{flushright}

December\qss 30,\qss 2014\linebreak

\vspace*{-1.5\bigskipamount}
April\qss 16,\qss 2016\qss (Minor edits and\qss ``Note bibiliographique'')

\vspace{\bigskipamount}
 
http:/\!/\hspace*{-0.07em}nikolaivivanov.com
\end{flushright}

\end{document}